\documentclass[12pt]{amsart}

\usepackage{amsmath,amssymb, amscd, stmaryrd, bm}
\usepackage{fullpage}
\usepackage{enumerate}
\usepackage{cite}
\usepackage{framed}
\usepackage[dvips,dvipdf]{graphicx}
\usepackage[usenames,dvipsnames]{color}
\usepackage{color}
\usepackage{tabularx}
\usepackage{array}
\usepackage{multirow}

\usepackage{amsfonts}
\usepackage{epstopdf}
\usepackage{enumerate}
\usepackage{cleveref}
\usepackage{subfig,caption,float}
\usepackage{mathrsfs}
\usepackage{booktabs}
\newcounter{example}[section]
\newenvironment{example}[1][]{\refstepcounter{example}\par\medskip
	\noindent \textbf{Example~\theexample. #1} \rmfamily}{\medskip}

\graphicspath{{figures/}}

\newcommand{\ncom}{\newcommand}
\ncom{\beqn}{\begin{eqnarray*}}
	\ncom{\eeqn}{\end{eqnarray*}}

\theoremstyle{definition}

\theoremstyle{definition}

\theoremstyle{remark}

\title{FAST ACCURATE APPROXIMATION OF CONVOLUTIONS WITH
	WEAKLY SINGULAR KERNEL AND ITS APPLICATIONS}

\author[A. K. Tiwari]{Awanish Kumar Tiwari} 
\address{Awanish Kumar Tiwari, Department of
	Mathematics, Indian Institute of Science Education and Research Bhopal, MP 462066}
\email{awanisht@iiserb.ac.in}

\author[A. Pandey]{Ambuj Pandey} 
\address{Ambuj Pandey, Department of Mathematics, Indian Institute of Science Education and Research Bhopal, MP 462066}
	\email{ambuj@iiserb.ac.in}
	
\author[J. Paul]{Jagabandhu Paul}
\address{Computational \& Mathematical Sciences, California Institute of Technology} \email{jpaul@caltech.edu}

\author[A. Anand]{Akash Anand} 
\address{Akash Anand, Department of
	Mathematics and Statistics, Indian Institute of Technology, Kanpur, UP 208016}
\email{akasha@iitk.ac.in}

\begin{document}
\maketitle{\title}
\begin{abstract}
	In this article, we present an $O(N \log N)$ rapidly convergent algorithm for the numerical approximation of the convolution integral with radially symmetric weakly singular kernels and  compactly supported densities.
	To achieve the reduced computational complexity, we utilize the Fast Fourier Transform (FFT)  on a uniform grid of size $N$ for approximating the convolution. To facilitate this and maintain the accuracy,  we primarily rely on a periodic Fourier extension of the density with a suitably large period depending on the support of the density.
	The rate of convergence of the method increases with increasing smoothness of the periodic extension and, in fact, approximations exhibit super-algebraic convergence when the extension is infinitely differentiable. Furthermore, when the density has jump discontinuities, we utilize a certain Fourier smoothing technique to accelerate the
	convergence  to achieve the quadratic rate in the overall approximation.   Finally, we apply the integration scheme for numerical solution of certain partial differential equations. Moreover, we apply the quadrature to obtain a fast and high-order Nyst\"om solver for the solution of the Lippmann-Schwinger integral equation. We validate the performance of the proposed scheme in terms of accuracy as well as computational efficiency through a variety of numerical experiments.	
\end{abstract}
\section{Introduction}

In this article,  we propose an $O(N \log N)$ rapidly convergent algorithm  for numerical approximation of the convolution operator of the form
\begin{equation}~\label{eq:conv}
(Au)(x) = \displaystyle \int_D g(x-y)u(y)\,dy
\end{equation}
with a weekly singular kernel \(g\) and a smooth density $u$ whose support is contained in $D$ which, without loss of generality, we assume to be the unit square, that is, $D= [0,1]\times[0,1]$. 
The convolution integral of the form \cref{eq:conv} arise in a wide range of application areas, for instance, 
fluid dynamics \cite{19347}, quantum physics \cite{bao_tang_zhang_2016}, and electromagnetism \cite{bruno2004efficient}, to name a few. In particular, utilizing the proposed numerical scheme, we construct fast solvers for two important problems. For the first set of problems, we obtain solutions to certain linear constant-coefficient partial differential equations that arise, for example, in electrostatics and fluid dynamics. For the second problem, we solve the Lippmann-Schwinger integral equation whose application areas include wave propagation and scattering \cite{colton1998inverse}, medical imaging \cite{10.2307/118292, colton1999linear, pham2020three}, and Radar imaging \cite{2019RS006946}. 

In the literature, several fast, high-order algorithms are available, especially in the context of Helmholtz and Laplace kernels. 
In particular, most convolutions where \(g(x)\) is the fundamental solution of the Laplace equation arise in the context of the Poisson problem \cite{ethridge2001new}.   
For a concise review on the more general problem of numerical solution of constant-coefficient partial differential equations, we refer the readers to  \cite{BAO201572, doi:10.1063/1.2335442,exl2016accurate}. On the other hand, most solvers involving the Helmholtz kernel are devoted to the numerical solution of the Lippmann-Schwinger integral equation 
\cite[page 216]{colton1998inverse}.
While we do not intend to review all such work here, some recent contributions in this direction can be found in \cite{ethridge2001new,beylkin2009fast,
	imbert2019integral,duan2009high,marin2014corrected,aguilar2004high,bruno2004efficient,andersson2005fast,hesford2010fast,DARVE2000195, 
	greengard_rokhlin_1997,250128} 
and references therein. 
To the best of our knowledge, none of the methods cited above are designed to handle
the more general class of weakly singular kernels.
A notable exception is the method proposed in 
\cite{vico2016fast} that can handle a large class of weakly singular kernels.   The key idea in \cite{vico2016fast} is to truncate the 
underlying kernel 
so that the value of convolution integral \cref{eq:conv} remains unchanged.  As a result, by the Paley-Wiener theorem, the Fourier transform of the whole integrand is infinitely smooth, and subsequently, an accurate approximation of \cref{eq:conv} is obtained by evaluating the inverse Fourier transform of smooth functions. While the method is easy to implement and converges fast, its applicability is constrained by the availability of an analytical expression of the Fourier transform of the underlying truncated singular kernels, which, in general, is not available.

%
%
%
%

The proposed numerical scheme achieves the desired $O(N \log N )$ computational complexity by employing the Fast Fourier Transform (FFT) based on equispaced Cartesian grids and is applicable to any weakly singular kernel. The method converges super-algebraically provided the density $u$ is infinitely differentiable in $D$ and all its derivatives vanish at the boundary of $D$. As we see in the next section, the proposed integration scheme requires computation of certain moments. For some special kernels, including those coming from the fundamental solution of the Laplace and the Helmholtz equation, analytic expressions are available for these moments. In such cases, our method reduces to a straightforward application of the two-dimensional FFT. However, more generally, when such exact expressions are not available, we can obtain them numerically to a desired level of accuracy. For example, we discuss one such approach to pre-compute their high-order numerical approximations. Even though this pre-computation is needed to be performed only once for each weakly singular kernel, the accompanying cost can be reduced dramatically by invoking a specific localization technique.
Moreover, in the case when $u$ is only piecewise smooth, we utilize a certain Fourier smoothing technique \cite{bruno2004efficient, pandey2020fourier} to enhance the linear rate of convergence to achieve approximations that converge quadratically.

%
%

The organization of the paper is as follows. We present the proposed numerical integration scheme to approximate the convolution integral \eqref{eq:conv} in \Cref{sec:descr-quadr}. We begin, in \Cref{sec:glob-quad}, with an introduction to the numerical method and follow it up,  in \Cref{sec:weight}, with a discussion on the computation of certain moments that are needed in the method. 
In \Cref{sec:local-beta-weights}, we discuss the localized integration scheme to reduce the pre-computation cost when an analytic expression for obtaining  these moments is not available.
Next, to improve the convergence rate for piecewise smooth densities, a Fourier smoothed integration scheme is discussed in \Cref{sec:FS}. Finally, we present a variety of numerical results to validate the performance and applicability of the proposed method in \Cref{sec:numerical-examples}. In particular, in \Cref{sec:convol}, we conduct extensive numerical experiments to computationally establish the expected convergence rates of our scheme. Moreover, we also include a test to confirm the $O(N\log N)$ computational complexity of the method. Further, in \Cref{seq:poison} and \Cref{sec:appl-volum-scatt}, we apply the integration scheme to solve  Poisson and Lippmann-Schwinger equation, respectively, in a variety of settings. Some concluding remarks along with some future directions are put forth in \Cref{sec:conclusions}.

\section{The numerical integration scheme}\label{sec:descr-quadr}
This section presents a detailed description of our numerical scheme for the approximation of \cref{eq:conv}. 

\subsection{Basic Idea}\label{sec:glob-quad}
Toward obtaining an approximation to \cref{eq:conv}, we begin by considering the $b$-periodic extension $u_e$ of $u$ by
\begin{equation}
\label{eq:fc_discrete}
u_e(x) = 
\begin{cases}
u(x), & x \in D \\
0, & x \in D_e \setminus D,
\end{cases}
\end{equation}
where $D_e = [0,b]^2$, 
and $b$ is an integer with $b \ge 1+\sqrt{2}$. If $B_a(x)$ denotes the ball of radius $a$ centered at $x$, then 
for all $a$ satisfying $\sqrt{2} \le a \le b-1$ and all $x \in D$, 
we have $D \subseteq B_a(x) \subseteq [-(b-1),b]^2$, so that
\begin{align} \label{eq:io}
(Au)(x) = \int_D g(x-y)u(y)dy = \int_D g(x-y)u_e(y)dy = \int_{B_a(x)} g(x-y)u_e(y)dy.
\end{align}
Using the $b$-periodicity of \(u_e\),
the last integral in \cref{eq:io} is expressed as
\begin{align}~\label{eq:A} 
(Au)(x) 
&= \sum_{k \in \mathbb{Z}^2} \widehat{g}(k) \widehat{u_e}(k) e^{2\pi i k\cdot x/b},
\end{align}
where
\begin{align}\label{eq:bk}
\widehat{g}(k) &= 
\displaystyle \int_{B_a(0)}
g(y) e^{-2\pi ik \cdot y/b} \,dy,\ \ \ \widehat{u_e}(k)
= \dfrac{1}{b^2}\displaystyle \int_{D_e} u_e(y)
e^{-2\pi ik\cdot y/b}dy.
\end{align}

In view of \cref{eq:A} and \cref{eq:bk}, we take the numerical integration scheme $A_n$ to be of the form 
\begin{align}\label{eq:An}
(A_n u)(x)
&= \sum_{k \in \mathbb{F}_n }
\widehat{g}(k) \widehat{u_{e,n}}(k)e^{2\pi i k\cdot x/b},
\end{align}
where 
\[
\mathbb{F}_n = 
\begin{cases}
\mathbb{Z}^2 \cap [-nb/2,nb/2)^2, & nb \text{ is even},  \nonumber\\
\mathbb{Z}^2 \cap [-(nb-1)/2,(nb-1)/2]^2, &  nb \text{ is odd},
\end{cases}
\]
and $\widehat{u_{e,n}}(k)$ is the trapezoidal approximation of $\widehat{u_{e}}(k)$ given by
\begin{align*}
\widehat{u_{e,n}}(k) = \frac{1}{(nb)^2} \sum_{j \in \mathbb{G}_n} u_j e^{-2\pi ik\cdot j/nb},
\end{align*}
where $\mathbb{G}_n = \mathbb{Z}^2 \cap ([0,n)\times[0,n))$ and $u_j = u(j/n)$ is the grid data on the regular mesh $\{ j/n, j \in \mathbb{G}_n\}$ on $D$.
Note that the evaluation of the coefficients $\widehat{u_{e,n}}$ can be accomplished in $O(N\log N)$, $N = (bn)^2$, operations using FFT.
The numerical scheme \cref{eq:An} can be rewritten in the form of a quadrature 
given by
\begin{align}
(A_nu)(x)
&=  \sum_{j \in \mathbb{G}_n} w_j^n(x) u(x_j)
\label{eq:quad}
\end{align}
with quadrature points $x_j = j/n$, and weights 
\begin{align} \label{eq:quad_wts}
w_j^n(x) = \frac{1}{(nb)^2} \sum_{k \in \mathbb{F}_n} \widehat{g}(k)    e^{2\pi ik\cdot \left(x - x_j\right)/b}, \ \ j \in \mathbb{G}_n.
\end{align}

\subsection{Moments}\label{sec:weight}

While \cref{eq:An} provides an $O(N\log N)$ scheme for obtaining the convolution on an equispaced grid of size $N = (bn)^2$, the quadrature  \cref{eq:quad} with the pre-computed weights can be used to compute $(A_nu)(x)$ for any given $x$. Use of either form requires that the moment $\widehat{g}(k)$ is available for all $k \in \mathbb{F}_n$. In particular, for a radially symmetric $g$, that is, $g(y)=g(|y|)$, the calculation of $\widehat{g}(k)$ simplifies as follows
\begin{align} \label{eq:ghat}
\widehat{g}(k) = \int_{B_a(0)} g(y)e^{-2\pi ik \cdot y/b}\,dy = 2\pi \int_0^a g(\rho) J_0(2\pi |k| \rho/b)\rho\,d\rho,
\end{align}
where the radial integral can be evaluated analytically in several cases of practical interest.
For instance, if $g(x) = -\frac{1}{2\pi}\log|x|$ (Laplacian kernel),  we have
\begin{align}
\widehat{g}(k) &
= 
\begin{cases}
\dfrac{a^2}{4}(1-2 \log a), & |k|=0,  \nonumber \\
- \dfrac{a\log a\ J_1(2\pi |k| a/b)}{2\pi |k|/b} + \dfrac{1-J_0(2\pi |k| a/b) }{(2\pi |k|/b)^2} , & |k|\ne 0.
\end{cases} 
\end{align} 
Similarly, for $g(x) = \frac{i}{4}H_0^{(1)}(\kappa |x|)$, a kernel of interest for scattering calculations in two dimensions, we have
\begin{align*} 
&\widehat{g}(k) = 2\pi \int_0^a H_0^{(1)}(\kappa \rho)J_0(2\pi |k| \rho/b)\rho\,d\rho =  \\
&
\begin{cases}
\frac{\pi}{2\kappa}\left(i a J_1(\kappa a)-a Y_1(\kappa a) + \frac{2}{\kappa \pi}\right), & |k|=0, \nonumber \\
\dfrac{i\pi}{2} \dfrac{\left(\frac{2\pi |k| a}{b}\right) J_1\big(\frac{2\pi |k|a}{b}\big)H_0^{(1)}(\kappa a)-\kappa a J_0\big(\frac{2\pi |k|a}{b}\big)H_1^{(1)}(\kappa a) - \frac{2i}{\pi}}{(2\pi |k|/b)^2 - \kappa^2}, & \kappa \ne \frac{2\pi |k|}{b},  |k|\ne 0,  \nonumber\\
\dfrac{i}{4} \pi a^2\left( J_0(\kappa a)H_0^{(1)}(\kappa a) + J_1(\kappa a)H_1^{(1)}(\kappa a) \right), &  \kappa =  \frac{2\pi |k|}{b}, |k|\ne 0.
\end{cases}
\end{align*}
For $g(x) = (2\pi |x|)^{-1}$, the values for $\widehat{g}(k) $ is given by the expression
\begin{align*}
\begin{cases} 
a, & |k|=0, \nonumber \\
a\left( J_0\big(\frac{2\pi |k|a}{b}\big) + \frac{\pi}{2}\left( J_1\big(\frac{2\pi |k|a}{b}\big)\boldsymbol{H}_0\big(\frac{2\pi |k| a}{b}\big) - J_0\big(\frac{2\pi |k|a}{b}\big)\boldsymbol{H}_1\big(\frac{2\pi |k| a}{b}\big)\right)\right), & |k|\ne 0,
\end{cases}
\end{align*}
where $\boldsymbol{H}_{\nu}(x)$ denotes the Struve function of order $\nu$.
The corresponding expression for $g(x) = \frac{1}{2\pi}|x|^{-1/2}$ reads
\begin{align*}
\begin{cases} 
\frac{2}{3} a^{3/2}, & |k|=0, \nonumber \\
\frac{\sqrt{2}}{ (2\pi |k|/b)^{3/2}}\frac{\Gamma(3/4)}{\Gamma(1/4)} - \\
\frac{a}{(2\pi |k|/b)^{1/2}}\left(\frac{1}{2}J_0\big(\frac{2\pi |k|a}{b}\big)S_{-1/2,-1}\big(\frac{2\pi |k| a}{b}\big) - J_1\big(\frac{2\pi |k|a}{b}\big)S_{1/2,0}\big(\frac{2\pi |k| a}{b}\big)\right), & |k|\ne 0,
\end{cases}
\end{align*}
where $S_{\mu,\nu}$ are Lommel functions (see \cite{gradshteyn2007}).
More generally, for $g(x) = \frac{1}{2\pi}|x|^{\gamma}, \gamma > -2$, the corresponding expression is given by
\begin{align*}
\begin{cases}
\frac{1}{2+\gamma} a^{2+\gamma}, & |k|=0, \nonumber \\
\frac{2^{1+\gamma}}{ (2\pi |k|/b)^{2+\gamma}}\frac{\Gamma(1+\gamma/2)}{\Gamma(-\gamma/2)} + \\
\frac{a}{(2\pi |k|/b)^{1+\gamma}}\left(\gamma J_0\big(\frac{2\pi |k|a}{b}\big)S_{\gamma,-1}\big(\frac{2\pi |k| a}{b}\big) + J_1(\frac{2\pi |k|a}{b})S_{1+\gamma,0}\big(\frac{2\pi |k| a}{b}\big)\right), & |k|\ne 0.
\end{cases}
\end{align*}
However, for a general weakly singular kernel, 
an analytic expression for $\widehat{g}(k)$ may not be readily available. In such cases where obtaining exact expression is not possible, we can numerically pre-compute them with high precision for their later use. The numerical computations of  $\widehat{g}(k)$ poses difficulties 
mainly due to the singularity present in the integrand at \(\rho=0\). We resolve this singularity
using a change of variable of the form \(\rho =\tau^p\) by choosing a sufficiently large integer $p>1$.
The integrand of resulting integral
\begin{align*}
\widehat{g}(k)
&= 2 \pi p \int\limits_{0}^{a^{1/p}}
g(\tau^p) J_0(2\pi |k|\tau^p/b) \tau^{2p-1}
\,d\tau,
\end{align*}
has a high degree of smoothness in $\tau$, and therefore, can be computed accurately using a high-order quadrature. For calculations in this paper, we used the Clenshaw-Curtis quadrature where, for integration domain $[0,\alpha]$, with $\alpha > 0$, the quadrature points $\{x^{cc}_j: j = 1,\ldots,n_{pc}\}$ and the weights $\{\omega^{cc}_j: j = 1,\ldots,n_{pc}\}$
are given by 
\begin{align} \label{eq:cc_pts}
x^{cc}_j = \frac{\alpha}{2}\left(1+ \cos\left(\left(j-\frac{1}{2}\right)\frac{\pi}{n_{pc}}\right)\right),
\end{align}
and 
\begin{align} \label{eq:cc_wts}
\omega^{cc}_j = \frac{\alpha}{n_{pc}} \sideset{}{'}\sum_{k=0}^{n_{pc}} \vartheta_k \cos\left(\left(j-\frac{1}{2}\right)\frac{k\pi}{n_{pc}}\right),
\end{align}
with
\begin{align*}
\vartheta_k = 
\begin{cases}
-2/(k^2-1), & \text{if k is even}, \\
0 , & \text{if k is odd},
\end{cases}
\end{align*}
where the primed sum denote that the first term in the summation is halved.

\subsection{Localization}
~\label{sec:local-beta-weights}
In many large calculations, especially when an analytic expression is not available for $\widehat{g}$, one can reduce the pre-computation time by localizing the singularity in the integral \cref{eq:bk} defining
\(\widehat{g}(k)\) to a small neighborhood of the singularity. 
Toward this, we start by splitting the integral operator
\((A u)(x)\) using a smooth windowing function, say $\eta_{\beta}$ for some $\beta > 0$, with the following properties -- (a) $\eta_{\beta} : [0,\infty) \to [0,1]$ is an infinitely differentiable function, (b) $\eta_{\beta}(|x|) = 0$ for all $x$ with $|x| \ge {\beta}$ and (c) $\eta_{\beta}(|x|) = 1$ for all $x$ with $|x| \le \beta /2$. Indeed, using $\eta_{\beta}$, the convolution integral in \cref{eq:conv} can be split into two integrals as
\begin{align}
\label{eq:io_split}
(A u)(x) 
&= \int_{D} g(x-y) \big(1-\eta_{\beta}(|x-y|)\big)u(y) \,dy 
+ \int_{B_{\beta}(x)} g(x-y)  \eta_{\beta}(|x-y|)u(y)\,dy.
\end{align}
The first convolution in \cref{eq:io_split} with the smooth kernel $g(x) \big(1-\eta_{\beta}(|x|)\big)$ and compactly supported density $u$ is approximated to high-order using  the classical trapezoidal rule. Obviously, when this is implemented using the FFT, the computational cost reduces to $O(N\log N)$.
The second convolution in \cref{eq:io_split} with the localized weakly singular kernel $g(x)  \eta_{\beta}(|x|)$ is computed using the numerical integration scheme given in \cref{eq:An}, where $\widehat{g}$ is appropriately replaced with  
\begin{align*}
\widehat{g}_{\eta_{\beta}}(k)
&= \int_{B_\beta(0)} g(y) \eta_{\beta}(|y|) e^{-2\pi ik\cdot y/b}\,dy
\end{align*}
which, for a radially symmetric $g$, reduces to
\begin{align*}\label{eq:-Bk2d}
\widehat{g}_{\eta_{\beta}}(k)
&= 2 \pi \displaystyle \int_{0}^{\beta}g(\rho) \eta_{\beta}(\rho) J_0(2\pi |k|\rho/b)\rho\,d\rho. 
\end{align*}
Again, we can pre-compute these localized weights numerically using a high-order quadrature after resolving the singularity through the change of variable \(\rho =\tau^p\) for sufficiently large integer $p>1$. Although we can again employ the quadrature with points and weights given by \cref{eq:cc_pts}-\cref{eq:cc_wts}, however, as the integrand vanishes to a high-order at both endpoints, the use of the classical trapezoidal rule is as accurate and computationally cheaper. Moreover, the relatively smaller integration domain reduces the pre-computation time even further while maintaining the desired accuracy.

%
%

\subsection{Improved rate of convergence for discontinuous density}
~\label{sec:FS}
While the integration scheme presented in \Cref{sec:glob-quad} converges to high-order for smooth and compactly supported density $u$, it converges with only first-order when $u$ is discontinuous in $D$.
In this section, we explore the Fourier smoothing technique introduced in \cite{bruno2004efficient} and further developed in \cite{pandey2020fourier,bruno2019fast}, to enhance the convergence rate  while maintaining $O(N \log N)$ computational cost.
%
Use of the resulting quadrature rule in the Nystr\"{o}m method for the solution of the  Lippmann-Schwinger equation, in turn,  
produces second-order convergent approximations when material discontinuities are present in the scattering geometry. 
This enhancement in the convergence rate is a notable improvement as most of the existing fast numerical techniques yield only the first-order convergent solution when we have discontinuity across the material interfaces. We know that the Fourier smoothing also alleviates numerical difficulties that arise due to non-smooth boundaries \cite{bruno2004efficient}, a property which remains true for the proposed scheme as well, as seen in \Cref{exp:9} in the next section.

\begin{figure}[h!] 
	\centering
	\subfloat[Surface view of $\psi$] {
		\includegraphics[scale =0.21]{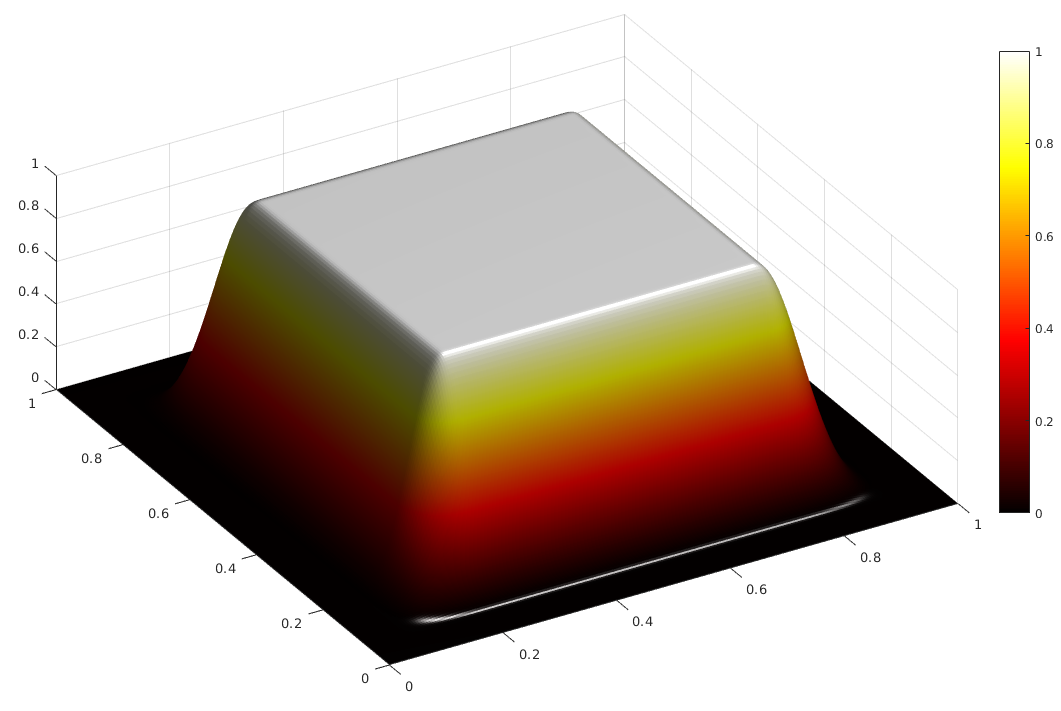}
	}\label{fig:eta1}
	\subfloat[1-d view of $\psi$] {
		\includegraphics[scale=0.36]{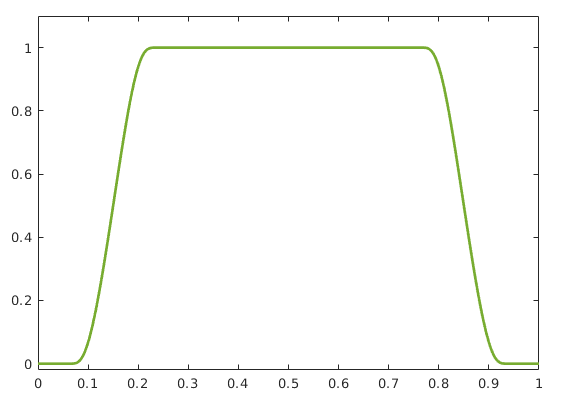}
	}\label{fig:eta2}
	\caption{Smooth vanishing window function $\psi$. }\label{fig:eta}
\end{figure}

More precisely, we consider the indicator function $\chi_{\Omega}(x)$ satisfying $\chi_{\Omega}(x)=1$ if $x \in \Omega=\{x \ | \ u(x)\neq 0\}\subsetneq D,$ and zero otherwise. Let $\tilde{u}$ be a smooth extension of $u$
such  that $\tilde{u}(x)=u(x)$ for all $x \in \Omega$. Now, using $\tilde{u}$ and $\chi_{\Omega}$, we recast the
integral operator~\cref{eq:conv} as
\begin{equation}\label{ls1}
(Au)(x) = \int_D g(x-y) \tilde{u}(y)\chi_{\Omega}(y)\psi(y)dy,
\end{equation}
where $\psi$ is a smooth window function such that $\psi(x) = 1$ for $x \in \Omega$ and $\psi(x) = 0$ for $x \in \partial D$.
For illustration, a specific instance of $\psi$ is depicted in Figure~\ref{fig:eta}. 
The key idea is to replace $\chi_{\Omega}$ in \cref{ls1} by its truncated Fourier series 
\begin{equation}  \label{tfs} 
\chi_{\Omega}^{n}(x) = \sum_{k \in \mathbb{F}_n } \widehat{\chi_{\Omega}^{n}}(k) e^{2\pi i k\cdot x}.
\end{equation}
Using the integration scheme $A_n$, we obtain the approximation to $Au$ as $A_n(\tilde{u}\chi^n_{\Omega}\psi)$, 
that converges with second-order provided 
Fourier coefficients $\widehat{\chi_{\Omega}^{n}}(k)$ are known  accurately either numerically or analytically. A partial mathematical justification for the second-order convergence by means of this Fourier smoothing technique can be found in \cite{bruno2005higher, hyde2005fast}.


\section{Numerical Examples}
~\label{sec:numerical-examples}
As we have mentioned in the introduction, this paper aims twofold: first, to design a fast high order integration scheme to approximate convolution integral \cref{eq:conv}, and second, its application to the solution of constant linear coefficient partial differential equations and Lippmann-Schwinger integral equation. First, in \Cref{sec:convol}, we illustrate our integration scheme's high order character and performance for various weakly singular kernels with different degrees of smoothness in the extended density $u_e$. Subsequently, in \Cref{seq:poison} and \Cref{sec:appl-volum-scatt}, we present simulation results for certain partial differential equations and the Lippmann-Schwinger equation.

The relative error reported in all the examples are computed as
\begin{equation}~\label{eq:rel-err}
\varepsilon_{n} = \frac{\max_{j \in \mathbb{G}_n}|(\mathrm{U}_{\mathrm{exact}})(j/n) - (\mathrm{U}_{\mathrm{approx}})(j/n)|}{\max_{j \in \mathbb{G}_n}|(\mathrm{U}_{\mathrm{exact}})(j/n)|},
\end{equation}
where $\mathrm{U}_{\mathrm{approx}}$ and $\mathrm{U}_{\mathrm{exact}}$ (either obtained analytically or computed by our scheme at a very fine grid) 
denote the approximate and exact  values of the quantity to be computed, respectively. In the formula \cref{eq:rel-err}, subscript $n$ stands for the number of grid points in each dimension. 
The numerical order of convergence (noc),  which relates to the rate at which the accuracy of the approximation improves as $n$ increases, is computed by the following formula:
%
\begin{equation}~\label{eq:ord-conv}
\textrm{noc}
= \log_{2}\left(\frac{\varepsilon_{n}}{\varepsilon_{2n}}\right).
\end{equation}  
For all experiments in this section, we have taken the extended computational domain $D_e = [0,b]^2$ with  $b=3$, the smallest integer greater than $1+\sqrt{2}$.
%

\subsection{Approximation of the convolution integral}~\label{sec:convol}

In this section, we demonstrate the performance (convergence rate and computational complexity) of the proposed integration scheme for approximation of the convolution integral \cref{eq:conv} through a variety of examples.
\begin{example}\label{exp:smooth}
	\emph{(Super-algebraic convergence)}		 
	In this example, we demonstrate the performance of our numerical integration scheme for infinitely smooth compactly supported density for several weakly singular kernels. 
	We obtain numerical approximations of \cref{eq:conv} for four different singular kernels, namely, $g(x) = |x|^\gamma,\, \gamma = -1/2, -1, -3/2$, and $g(x) = \log(|x|)$ with density
	\begin{equation}\label{gaus}
	u(y) = (1/2\pi \sigma^2)\exp(-(|y|^2/2\sigma^2)),
	\end{equation}
	where $\sigma = 0.05$. The results are presented in \cref{table:globconv_gauss} which  clearly corroborates the super algebraic convergence of the proposed algorithm for infinitely smooth and compactly supported densities.  		 		
	
	\begin{table}[t!]
		\begin{center}
				\begin{tabular}{ c | c | c |c |c |c |c |c |c }
					\hline
					& \multicolumn{2}{l|}{$g(x) = |x|^{-1/2}$} & \multicolumn{2}{l|}{$g(x) = |x|^{-1}$} & \multicolumn{2}{l|}{$g(x) = |x|^{-3/2}$}& \multicolumn{2}{l}{$g(x) = \log(|x|)$}\\ \hline
					n & $\varepsilon_{n}$ & noc & $\varepsilon_{n}$ & noc & $\varepsilon_{n}$ & noc& $\varepsilon_{n}$ & noc\\ \hline
					$2^2$ & $ 1.7\times10^{0}$ & $-$ & $ 8.9\times10^{-1}$ & $-$ 
					& $ 3.9\times10^{-1}$ & $-$& 
					$ 2.4\times10^{0}$ & $-$\\ 
					
					$2^3$ & $ 2.3\times10^{-1}$ & $3.1$ & $ 1.8\times10^{-1}$ & $2.6$& 
					$ 1.1\times10^{-1}$ & $2.0$& 
					$ 2.3\times10^{-1}$ & $3.6$\\ 
					
					$2^4$ & $ 2.7\times10^{-3}$ & $ 6.4$ & $ 1.7\times10^{-3}$ & $ 6.7$& $ 
					1.7\times10^{-3}$ & $ 6.0$& 
					$ 1.3\times10^{-3}$ & $ 7.5$\\ 
					
					$2^5$ & $1.6\times10^{-7}$ & $ 14.0$ & $1.1\times10^{-8}$ & $17.3$& 
					$1.5\times10^{-8}$ & $ 16.8$& 
					$3.8\times10^{-9}$ & $18.4$\\ 
					
					$2^6$ & $5.3\times10^{-15}$ & $24.9$ & $2.9\times10^{-16}$ & $25.1$ & 
					$6.6\times10^{-16}$ & $24.4$& 
					$2.5\times10^{-15}$ & $20.5$\\ \hline
				\end{tabular}
			\caption{Convergence study: Super algebraic convergence for a set of different weakly singular kernels with infinitely smooth compactly supported function  $u(y) = (1/2\pi \sigma^2)\exp(-|y|^2/2\sigma^2)$ with $\sigma = 0.05$.}\label{table:globconv_gauss}
		\end{center}
	\end{table}
	To numerically validate the computational complexity of our algorithm,  in \Cref{fig:complexity}, we have plotted the computational time (required for obtaining $(A_nu)(x_j),$ for all $j \in \mathbb{G}_n$) against the total number of grid points $N=(bn)^2$ for $n=24 \ell, \, \ell=1,\ldots,21$. For this simulation, we have taken $g(x) = \log|x|$ with the Gaussian density $u$ defined in \cref{gaus}.  
	The plot shows that the computational time indeed behaves as $C N\log N$, where $C = 1.2 \times 10^{-7}$,
	which confirms the expected $O(N\log N)$ complexity of the proposed integration scheme. 
	
	%
	\begin{figure}[h!] 
		\centering
		{
			\includegraphics[scale=0.4]{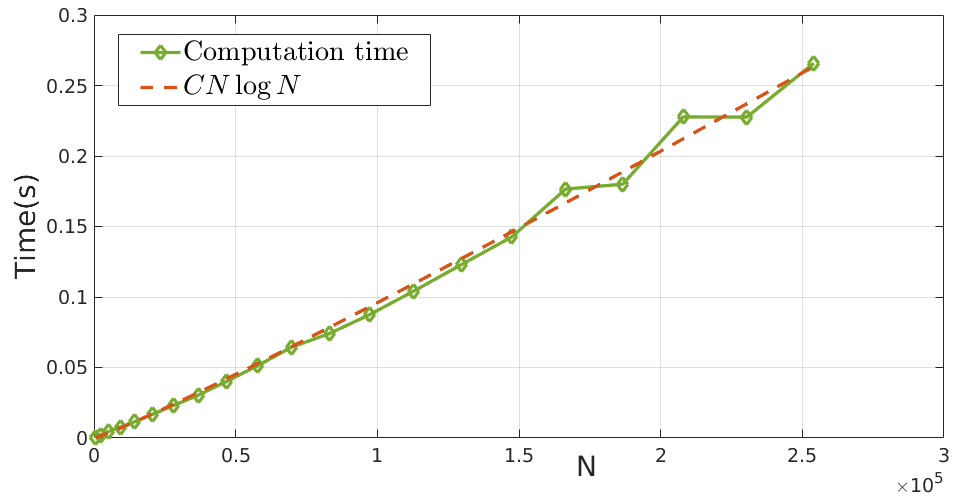}
		}
		\caption{The Computation time for the integration scheme against  $N=(bn)^2$. For comparison, a reference curve $C (N\log N)$ with  $C = 1.2\times10^{-7}$ is also shown.}\label{fig:complexity}
	\end{figure}
\end{example}

\begin{example}\label{exp:2}
	\emph{(The rate of convergence)}
	The error in our integration scheme depends on the asymptotic behaviour of  $|\widehat{g}(k)|$ and $|\widehat{u_e}(k)|$ (given in \cref{eq:bk}). To explore these dependence,
	we include two different sets of experiments.
	In the first group of simulations, numerical approximations of \cref{eq:conv} are obtained for a fixed density $u(x) = (x_1(1-x_1)x_2(1-x_2))^3$ with the same four kernels used in  \Cref{exp:smooth}. The corresponding convergence results are displayed in \Cref{table:globconv1}. 
	
	In our second set of experiments, we take $g(x) = |x|^{-1}$ and compute the approximations of \cref{eq:conv} for four different densities  $u(x) = (x_1(1-x_1)x_2(1-x_2))^m, \, m = 1, 2, 3, 4$ and the corresponding results are reported in  \cref{table:globconv2}. We see that, with increasing smoothness of $u_e$ as $m$ increases, the corresponding rate of convergence improves proportionally. 
	
	
	
	\begin{table}[t!]
		\begin{center}
			\begin{tabular}{ c | c | c |c |c |c |c |c |c }
				\hline
				& \multicolumn{2}{l|}{$g(x) = |x|^{-1/2}$} & \multicolumn{2}{l|}{$g(x) = |x|^{-1}$} & \multicolumn{2}{l|}{$g(x) = |x|^{-3/2}$}& \multicolumn{2}{l}{$g(x) = \log(|x|)$}\\ \hline
				n & $\varepsilon_{n}$ & noc & $\varepsilon_{n}$ & noc & $\varepsilon_{n}$ & noc& $\varepsilon_{n}$ & noc\\ \hline
				$2^2$ & $ 2.2\times10^{-2}$ & $-$ & $ 6.8\times10^{-3}$ & $-$ 
				& $ 4.9\times10^{-3}$ & $-$& $ 8.2\times10^{-3}$ & $-$\\ 
				
				$2^3$ & $ 8.0\times10^{-4}$ & $4.8$ & $ 6.1\times10^{-4}$ & $3.5$& $ 4.9\times10^{-4}$ & $3.3$& $ 4.8\times10^{-4}$ & $4.1$\\ 
				$2^4$ & $ 4.7\times10^{-5}$ & $ 4.1$ & $ 4.7\times10^{-5}$ & $ 3.7$& $ 4.8\times10^{-5}$ & $ 3.3$& $ 3.1\times10^{-5}$ & $ 4.0$\\ 
				$2^5$ & $3.1\times10^{-6}$ & $ 4.0$ & $3.4\times10^{-6}$ & $3.8$& $4.6\times10^{-6}$ & $ 3.4$& $1.9\times10^{-6}$ & 
				$4.0$\\ 
				$2^6$ & $1.9\times10^{-7}$ & $4.0$ & $2.4\times10^{-7}$ & 
				$3.8$ & $4.2\times10^{-7}$ & $3.4$& $1.2\times10^{-7}$ & 
				$4.0$\\ 
				$2^7$ & $1.2\times10^{-8}$ & $4.0$ & $1.7\times10^{-8}$ & $3.8$
				& $3.8\times10^{-8}$ & $3.5$& $7.4\times10^{-9}$ & $4.0$\\ 
				$2^8$ & $7.6\times10^{-10}$ & $4.0$ & $1.2\times10^{-9}$ & $3.9$ & $3.4\times10^{-9}$ & $3.5$& $4.6\times10^{-10}$ & $4.0$\\ \hline
			\end{tabular}
			\caption{An experiment to study the effect of $g$ on the convergence rate: Convolutions with different kernels and fixed density $u(x) = (x_1(1-x_1)x_2(1-x_2))^3$.
			}\label{table:globconv1}
		\end{center}
	\end{table}
	
	
	\begin{table}[t!]
		\begin{center}
			\begin{tabular}{ c | c | c |c |c |c |c |c |c }
				\hline
				& \multicolumn{2}{l|}{$m=1$}& \multicolumn{2}{l|}{$m=2$} & \multicolumn{2}{l|}{$m=3$}& \multicolumn{2}{l}{$m=4$} \\ \hline
				$n$ & $\varepsilon_{n}$ & noc& $\varepsilon_{n}$ & noc& $\varepsilon_{n}$ & noc& $\varepsilon_{n}$ & noc \\ \hline
				$2^3$ & $ 1.4\times10^{-2}$ & $-$ & $ 2.6\times10^{-4}$ & $-$& $ 6.1\times10^{-4}$ & $-$& $ 1.0\times10^{-4}$ & $-$\\ 
				
				$2^4$ & $ 4.2\times10^{-3}$ & $1.8$ & $ 4.2\times10^{-5}$ & $2.6$& $ 4.7\times10^{-5}$ & $3.7$& $2.0\times10^{-6}$ & 
				$5.7$\\ 
				$2^5$ & $ 1.2\times10^{-3}$ & $ 1.8$ & $ 5.9\times10^{-6}$ & $ 2.8$& $ 3.4\times10^{-6}$ & $ 3.8$& $ 7.6\times10^{-8}$ & $ 4.7$\\ 
				$2^6$ & $3.4\times10^{-4}$ & $ 1.8$ & $7.9\times10^{-7}$ & $ 2.9$& $2.4\times10^{-7}$ & $ 3.8$& $2.6\times10^{-9}$ & 
				$4.9$\\ 
				$2^7$ & $9.4\times10^{-5}$ & $1.9$ & $1.0\times10^{-7}$ & $3.0$& $1.7\times10^{-8}$ & $3.8$& $8.5\times10^{-11}$ & $4.9$ 
				\\ 
				$2^8$ & $2.6\times10^{-5}$ & $1.9$ & $1.3\times10^{-8}$ & $3.0$& $1.2\times10^{-9}$ & $3.9$& $2.7\times10^{-12}$ & $5.0$
				\\ 
				$2^9$ & $7.1\times10^{-6}$ & $1.9$& $1.7\times10^{-9}$ & $3.0$ & $7.9\times10^{-11}$ & $3.9$& $8.7\times10^{-14}$ & $5.0$\\ 
				\hline
			\end{tabular}
			\caption{An experiment to study the impact of smoothness of $u_e$ on the convergence rate: Convolution with fixed kernel $g(x) = |x|^{-1}$ and densities $u(x) = (x_1(1-x_1)x_2(1-x_2))^m, \, m=1,2,3,4$.}\label{table:globconv2}
		\end{center}
	\end{table}	
	
\end{example}
\begin{example}\emph{(Performance comparison of localized and non-localized schemes)}
	As discussed in \Cref{sec:local-beta-weights}, in absence of any analytical formula for obtaining the moments $|\widehat{g}(k)|$, they can be precomputed rapidly to a high degree
	of accuracy
	using the localization technique.
	In this example, to compare the computational times of local and non-local schemes, we compute the  
	convolution integral \cref{eq:conv} with  $g(x) = \log|x|$ and $u(x) = (x_1(1-x_1)x_2(1-x_2))^3$ using both methods for several values of $n$. We choose an optimal value for $n_{pc}$ in the pre-computation of localized as well as non-localized moments in such a way that the respective approximations have comparable accuracy level as the one obtained using the analytic expression for $|\widehat{g}(k)|$ given in \Cref{sec:weight}. 
	The respective computational times are displayed in \Cref{fig:comparison}(a) and the corresponding error are shown in \Cref{fig:comparison}(b). 
	We see that, due to the reduced pre-computation cost for the localized  integration, there is a significant improvement in the overall computation time as compared to the non-local scheme. 
	
	Next, to confirm that there is no loss of accuracy in the local
	scheme compared to the one without localization, we repeat the experiments done in \cref{exp:2} with localized moments obtained numerically. 
	The corresponding errors are reported in \cref{table:Conv1_loc} and \cref{table:conv_loc}. We see that the accuracy levels and convergence behaviour are nearly identical to those seen in \cref{table:globconv1} and \cref{table:globconv2}, respectively. 
	

	\begin{figure}[h!] 
		\centering
		\subfloat[Computational time for localized and non-localized schemes.] {
			\includegraphics[scale=0.37]{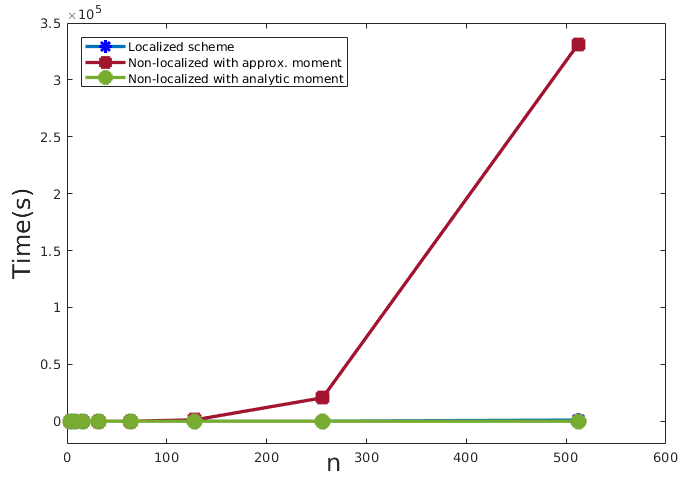}
		}\label{fig:Inttime}
		\subfloat[A log-log plot showing the errors in both approaches.] {
			\includegraphics[scale=0.385]{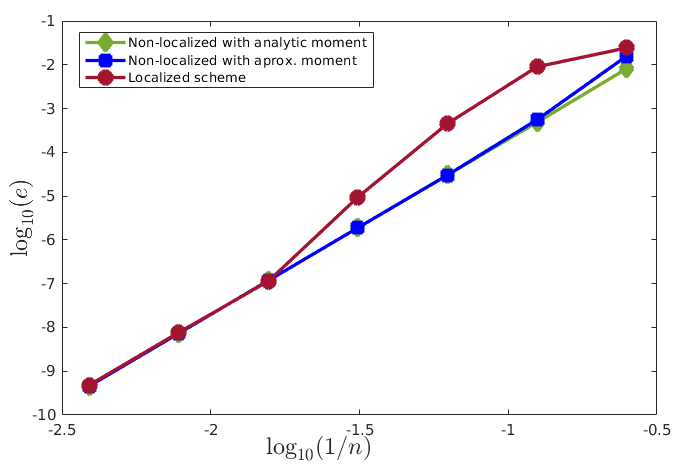}
		}\label{fig:err}
		\caption{A comparison of total computational time for localized and non-localized schemes for computing the convolution with $g(x) = \log(|x|)$ and $u(x) = (x_1(1-x_1)x_2(1-x_2))^3$: 
			The experiment indicate a significant improvement in the computational time for the localized scheme while maintaining the same degree of accuracy in the approximation.
		}\label{fig:comparison}
	\end{figure}
	
\end{example}

\begin{table}[b!]
	\begin{center}
		\begin{tabular}{ c | c | c |c |c |c |c |c |c }
			\hline
			& \multicolumn{2}{l|}{$g(x) = |x|^{-1/2}$} & \multicolumn{2}{l|}{$g(x) = |x|^{-1}$} & \multicolumn{2}{l|}{$g(x) = |x|^{-3/2}$}& \multicolumn{2}{l}{$g(x) = \log(|x|)$}\\ \hline
			$n$ & $\varepsilon_{n}$ & noc & $\varepsilon_{n}$ & noc & $\varepsilon_{n}$ & noc& $\varepsilon_{n}$ & noc\\ \hline
			$2^2$ & $ 7.6\times10^{-2}$ & $-$ & 
			$ 3.3\times10^{-2}$ & $-$ 
			& $ 5.2\times10^{-3}$ & $-$& 
			$ 2.4\times10^{-2}$ & $-$\\ 
			
			$2^3$ & $ 4.6\times10^{-3}$ & $4.0$ & 
			$ 4.1\times10^{-3}$ & $3.0$& 
			$ 4.9\times10^{-4}$ & $3.1$& 
			$ 9.0\times10^{-3}$ & $1.5$\\ 
			
			$2^4$ & $ 8.8\times10^{-5}$ & $ 5.7$ 
			& $ 2.0\times10^{-4}$ & $ 4.3$& 
			$ 4.9\times10^{-5}$ & $ 3.6$& 
			$ 4.5\times10^{-4}$ & $ 4.3$\\ 
			
			$2^5$ & $3.5\times10^{-6}$ & $ 4.7$ 
			& $5.4\times10^{-6}$ & $5.2$& 
			$5.1\times10^{-6}$ & $ 3.2$& 
			$9.5\times10^{-6}$ & $5.6$\\ 
			
			$2^6$ & $2.3\times10^{-7}$ & $3.9$ & 
			$3.4\times10^{-7}$ & $4.0$ 
			& $4.8\times10^{-7}$ & $3.4$& 
			$1.1\times10^{-7}$ & $6.4$\\ 
			
			$2^7$ & $1.4\times10^{-8}$ & $4.0$ & 
			$2.7\times10^{-8}$ & $3.7$& 
			$4.5\times10^{-8}$ & $3.4$& 
			$7.7\times10^{-9}$ & $3.9$\\ 
			
			$2^8$ & $9.0\times10^{-10}$ & $4.0$ & 
			$1.9\times10^{-9}$ & $3.8$ &
			$4.0\times10^{-9}$ & $3.5$& 
			$4.8\times10^{-10}$ & $4.0$\\ \hline
		\end{tabular}
		\caption{ An experiment to study the impact of $g$ on the convergence rate of the localized scheme: Convolution with various kernels and fixed density $u(x) = (x_1(1-x_1)x_2(1-x_2))^3$. }\label{table:Conv1_loc}
	\end{center}
\end{table}

\begin{table}[t!]
	\begin{center}
		\begin{tabular}{ c | c | c |c |c |c |c |c |c }
			\hline
			& \multicolumn{2}{l|}{$m=1$}& \multicolumn{2}{l|}{$m=2$} & \multicolumn{2}{l|}{$m=3$}& \multicolumn{2}{l}{$m=4$} \\ \hline
			$n$ & $\varepsilon_{n}$ & noc& $\varepsilon_{n}$ & noc& $\varepsilon_{n}$ & noc& $\varepsilon_{n}$ & noc \\ \hline
			$2^2$ & $ 6.0\times10^{-2}$ & $-$ & 
			$ 2.0\times10^{-2}$ & $-$& 
			$ 3.3\times10^{-2}$ & $-$& 
			$ 9.1\times10^{-3}$ & $-$\\ 
			$2^3$ & $ 1.5\times10^{-2}$ & $2.0$ 
			& $ 4.6\times10^{-3}$ & $2.1$& 
			$ 4.1\times10^{-3}$ & $3.0$& 
			$1.1\times10^{-3}$ & $3.0$\\ 
			$2^4$ & $ 5.1\times10^{-3}$ & $ 1.6$ & 
			$ 2.4\times10^{-4}$ & $ 4.3$& 
			$ 2.0\times10^{-4}$ & $ 4.3$& 
			$ 1.7\times10^{-5}$ & $ 6.0$\\ 
			$2^5$ & $1.6\times10^{-3}$ & $ 1.7$ 
			& $1.4\times10^{-5}$ & $ 4.1$& 
			$5.4\times10^{-6}$ & $ 5.2$& 
			$2.6\times10^{-7}$ & $6.1$\\ 
			$2^6$ & $4.7\times10^{-4}$ & $1.7$ & 
			$1.5\times10^{-6}$ & $3.2$& 
			$3.4\times10^{-7}$ & $4.0$& 
			$3.5\times10^{-9}$ & $6.3$ \\ 
			$2^7$ & $1.4\times10^{-4}$ & $1.8$ & 
			$1.9\times10^{-7}$ & $2.9$& 
			$2.7\times10^{-8}$ & $3.7$& 
			$1.2\times10^{-10}$ & $4.8$\\ 
			$2^8$ & $3.9\times10^{-5}$ & $1.8$& 
			$2.4\times10^{-8}$ & $3.0$ & 
			$1.9\times10^{-9}$ & $3.8$& 
			$4.0\times10^{-12}$ & $5.0$\\ 
			\hline
		\end{tabular}
		\caption{ An experiment to study the impact of smoothness of $u_e$ on the convergence rate of the localized scheme: Convolution with $g(x) = |x|^{-1}$ and $u(x) = (x_1(1-x_1)x_2(1-x_2))^m, \, m=1,2,3,4$.}\label{table:conv_loc}
	\end{center}
\end{table}	

\begin{example} \emph{(A kernel that is not radially symmetric)}
	As the final example of this section, we take the kernel 
	\[
	g(x) = -\dfrac{1}{2\pi}\dfrac{x_1}{x_1^2+x_2^2}
	\]
	to demonstrate that the proposed method is not limited to only radially symmetric kernels. We  compute the approximations of \cref{eq:conv} for four different densities. As the first experiment, we take $$u(x) = 4 \alpha(\alpha |x-c|^2 -1) \exp(-\alpha|x-c|^2)$$ with $\alpha = 250$ and $c = (1/2,1/2)$ and compare the numerical approximations with the exact value
	\begin{align*}
	(Au)(x) =	-2 \alpha (x_1-c_1) \exp(-\alpha|x-c|^2).
	\end{align*}
	The results reported in the two leftmost columns labeled ``Gaussian" in \Cref{table:genkern} confirms the exponential convergence.
	Next, we take the density $$u(x) = \Delta \left( (x_1(1-x_1)x_2(1-x_2))^m \right),$$ for $m = 4, 5, 6,$ and compare the approximations with the exact convolution
	\begin{align*}
	(Au)(x) = -m (2x_1-1)(x_1(1-x_1))^{m-1}(x_2(1-x_2))^m, \, m=4,5,6.
	\end{align*}
	The errors and corresponding numerical order of convergence are reported in \Cref{table:genkern}. 

	\begin{table}[b!]
		\begin{center}
				\begin{tabular}{ c | c | c |c |c |c |c |c |c }
					\hline
					& \multicolumn{2}{l|}{Gaussian}& \multicolumn{2}{l|}{$m=4$} & \multicolumn{2}{l|}{$m=5$}& \multicolumn{2}{l}{$m=6$} \\ \hline
					$n$ & $\varepsilon_{n}$ & noc& $\varepsilon_{n}$ & noc& $\varepsilon_{n}$ & noc& $\varepsilon_{n}$ & noc \\ \hline
					$2^2$ & $-$ & $-$ & 
					$ 8.1\times10^{-2}$ & $-$& 
					$ 3.4\times10^{-1}$ & $-$& 
					$ 5.3\times10^{-1}$ & $-$\\ 
					$2^3$ & $ 2.0\times10^{+1}$ & $-$ 
					& $ 2.4\times10^{-2}$ & $1.8$& 
					$ 1.4\times10^{-2}$ & $4.3$& 
					$3.6\times10^{-3}$ & $6.6$\\ 
					$2^4$ & $ 1.6\times10^{-1}$ & $ 3.8$ & 
					$ 3.2\times10^{-3}$ & $ 2.7$& 
					$ 7.1\times10^{-4}$ & $ 4.2$& 
					$ 2.7\times10^{-4}$ & $ 3.7$\\ 
					$2^5$ & $6.2\times10^{-5}$ & $ 11.3$ 
					& $4.3\times10^{-4}$ & $ 2.9$& 
					$4.3\times10^{-5}$ & $ 4.0$& 
					$1.0\times10^{-5}$ & $4.7$\\ 
					$2^6$ & $2.7\times10^{-16}$ & $37.6$ & 
					$5.5\times10^{-5}$ & $3.0$& 
					$2.6\times10^{-6}$ & $4.0$& 
					$3.3\times10^{-7}$ & $4.9$ \\ 
					$2^7$ & $-$ & $-$ & 
					$7.0\times10^{-6}$ & $3.0$& 
					$1.6\times10^{-7}$ & $4.0$& 
					$1.1\times10^{-8}$ & $5.0$\\ 
					$2^8$ & $-$ & $-$& 
					$8.8\times10^{-7}$ & $3.0$ & 
					$1.0\times10^{-8}$ & $4.0$& 
					$3.3\times10^{-10}$ & $5.0$\\ 
					\hline
				\end{tabular}
			\caption{A convergence study for a non-radial weakly singular kernel: The convolution integral with $g(x) = -x_1/(2\pi(x_1^2+x_2^2))$ and the ``Gaussian" $u(x) = 4 \alpha(\alpha |x-c|^2 -1) \exp(-\alpha|x-c|^2)$ with $\alpha = 250$ and $c = (0.5,0.5)$, and $u(x) = \Delta \left( (x_1(1-x_1)x_2(1-x_2))^m\right), \, m = 4, 5, 6$.
			}\label{table:genkern}
		\end{center}
	\end{table}	
	
\end{example}	

\subsection{Numerical solution of partial differential equations}\label{seq:poison}
It is known that if $g$ is a fundamental solution of the linear constant-coefficient partial differential operator $L$, then 
\[
v(x) = \int_{\mathbb{R}^2} g(x-y)f(y) \,dy, \quad x \in \mathbb{R}^2,
\]
solves $Lv = f$ in $\mathbb{R}^2$ provided $f \in L^1(\mathbb{R}^2)$. Some important examples include $L = -\Delta, \, L = -\Delta + \kappa^2I,$ and $L = \Delta + \kappa^2I$, where $I$ is the identity operator. Thus, if the source distribution $f$ is compactly supported, say in $D$, then the proposed integration scheme can be utilized to compute the solution to $Lv = f$.  

%

\begin{example}\label{exp:posn}
	We consider the Poisson equation
	\begin{equation}\label{eq:poison}
	-\Delta u(x) = f(x),
	\end{equation} 
	where the corresponding fundamental solution is given by $g(x) = -\frac{1}{2\pi}\log |x |$.
	For
	\[
	f(x) = \sum_{i=1}^{\ell} (4 \alpha^2 |x-c_i|^2 -4 \alpha) \exp(-\alpha|x-c_i|^2),
	\]
	the exact solution is given by
	\[
	u(x) = \sum_{i=1}^{\ell} \exp(-\alpha |x-c_i|^2).
	\]
	We compute the solution for $\alpha = 250$ and
	$\ell = 3$, where $c_1 = (0.6,0.6)$, $c_2 = (0.5,0.5)$,
	and $c_3 = (0.35,0.6)$. The results are given in \cref{table:poison}, where we clearly see a super-algebraic convergence. A similar convergence study is also included in  \cref{table:poison} for $\alpha = 950$ and  $\ell = 10$  with $c_1 = (0.6,0.6)$, $c_2 = (0.5,0.5)$,
	$c_3 = (0.35,0.6)$, $c_4 = (0.6,0.8)$, $c_5 = (0.8,0.8)$, $c_6 = (0.25,0.5)$, $c_7 = (0.75,0.5)$, $c_8 = (0.25,0.25)$, $c_9 = (0.5,0.25)$, and $c_{10} = (0.75,0.25)$.  
	The corresponding solutions are shown in \cref{fig:poison_sol}.
	\begin{figure}[tbh!] 
		\centering
		\subfloat[$\ell =3$] {
			\includegraphics[scale=0.31]{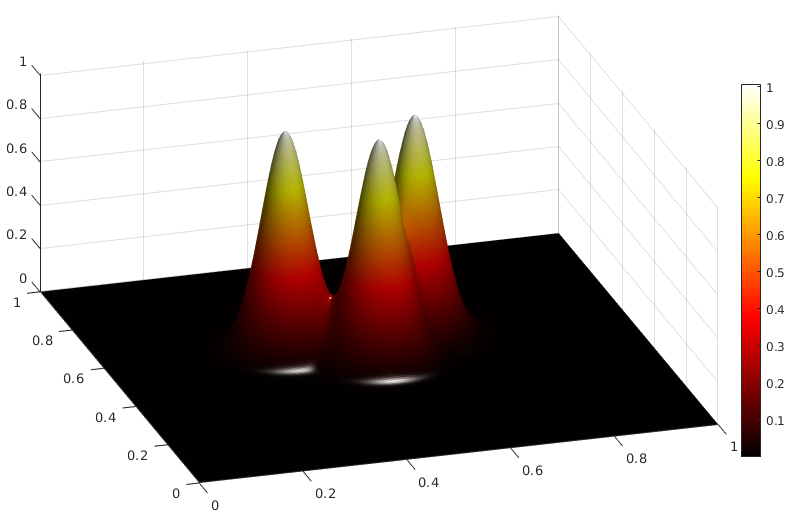}
		}\label{fig:poison1}
		\subfloat[$\ell =10$] {
			\includegraphics[scale=0.36]{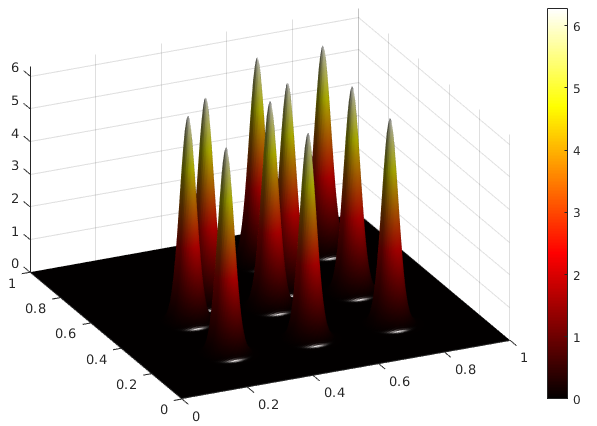}
		}\label{fig:poison2}
		\caption{Solution to the Poisson equation}\label{fig:poison_sol}
	\end{figure}
	\begin{table}[t]
		\begin{center}
			\begin{tabular}{ c | c | c  }
				\hline
				\multicolumn{3}{l}{$\ell = 3$}\\ \hline
				$ n$ & $\varepsilon_{n}$ & noc \\ \hline
				$ 2^2$ & $ 2.2\times10^{01}$ & $ -$ \\ 
				$ 2^3$ & $ 7.1\times10^{00}$& $ 1.6$ \\ 
				$ 2^4$ & $ 5.6\times10^{-02}$& $ 7.0$ \\ 
				$ 2^5$ & $ 2.8\times10^{-06}$ & $14.3$ \\ 
				$ 2^6$ & $ 9.7\times10^{-14}$ & $24.8$  \\ 
				\hline
			\end{tabular}
			\quad \quad
			\begin{tabular}{ c | c | c }
				\hline
				\multicolumn{3}{l}{$\ell = 10$}\\ \hline
				n & $\varepsilon_{n}$ & noc \\ \hline
				$2^3$ & $ 8.59\times10^{01}$ & $--$ \\ 
				$2^4$ & $ 1.47\times10^{01}$ & $2.6$ \\ 
				$2^5$ & $6.3\times10^{-02}$ & $ 7.9$ \\ 
				$2^6$ & $1.4\times10^{-06}$ & $15.5$ \\ 
				$2^7$ & $2.9\times10^{-15}$ & $28.5$ \\ \hline
			\end{tabular}
			\caption{A convergence study for the Poisson problem in  \Cref{exp:posn}.
			}\label{table:poison}
		\end{center}
	\end{table}
	
	We next consider $\Delta u = f$ with a discontinuous right hand side. More precisely, we take 
	\begin{equation}\label{eq:disc}
	f(x) = \begin{cases}
	1, & x \in [0.3,0.7]^2\\
	0, & x\in D\setminus [0.3,0.7]^2.
	\end{cases}
	\end{equation}	
	In this case, we compare the solutions obtained using the integration scheme with and without Fourier smoothing. 
	The results presented in \cref{table:conv_FS}  demonstrate a  significant improvement in the accuracy and enhancement in convergence rate for the Fourier smoothed scheme. For instance, for $n = 256$, 
	the result obtained from the Fourier smoothed version of the scheme is more than three thousand times accurate in comparison to that obtained without Fourier smoothing.
	\begin{table}[t!]
		\begin{center}
			\begin{tabular}{ c | c | c| c |c }
				\hline
				& \multicolumn{2}{l|}{WFS}& \multicolumn{2}{l}{FS}  \\ \hline
				$n$ & $\varepsilon_{n}$ & noc & $\varepsilon_{n}$& noc  \\ \hline
				$2^2$ & $ 3.6\times10^{-1}$ & $-$ & 
				$ 7.9\times10^{-2}$ & $ -$ \\ 
				$2^3$ & $ 1.9\times10^{-1}$ & $1.17$ 
				& $ 2.8\times10^{-3}$& $ 4.79$ \\ 
				$2^4$ & $ 7.2\times10^{-2}$ & $ 1.33$ & 
				$ 5.3\times10^{-4}$& $ 2.42$ \\ 
				$2^5$ & $3.7\times10^{-2}$ & $ 0.92$ 
				& $ 8.9\times10^{-5}$ & $2.57$ \\ 
				$2^6$ & $2.0\times10^{-2}$ & $0.91$ & 
				$ 2.1\times10^{-5}$ & $2.07$  \\ 
				$2^7$ & $9.9\times10^{-3}$ & $1.02$ & 
				$ 5.2\times10^{-6}$ & $2.04$ \\ 
				$2^8$ & $4.9\times10^{-3}$ & $1.02$& 
				$ 1.3\times10^{-6}$ &$2.02$ \\ 
				\hline
			\end{tabular}
			\caption{ A comparison between the Fourier smoothed (FS), and without Fourier smoothed (WFS) versions of the scheme for the discontinuous density given in \cref{eq:disc}. 
			}\label{table:conv_FS}
		\end{center}
	\end{table}		
\end{example}

\begin{example}
	Next, we consider the 
	equation 
	\begin{equation}\label{eq:yukawa}
	-\Delta u(x) + \kappa^2 u(x) = f(x),
	\end{equation} 
	with the fundamental solution 
	$g(x) =\frac{1}{2\pi}K_0(\kappa |x|)$, where $K_0$ is the modified Bessel function of the second kind of order $0$. For 
	\begin{equation}\label{eq:yukaw}
	f(x) = \left( \frac{-4|x|^2+4\delta^2}{\delta^4} +\kappa^2 \right) \exp\left(-\frac{1}{\delta^2}|x|^2\right),
	\end{equation}
	the exact solution is given by
	\[
	u(x) = \exp\left(-\frac{1}{\delta^2}|x|^2\right).
	\]
	The corresponding computational results for  $\delta = 0.08$ are shown in \cref{table:yukawa}, where we again observe the super algebraic convergence.
	
	\begin{table}[t]
		\begin{center}
			\begin{tabular}{ c | c | c  }
				\hline
				\multicolumn{3}{l}{$\kappa = 1$}\\ \hline
				$ n$ & $\varepsilon_{n}$ & noc \\ \hline
				$ 2^2$ & $ 1.5\times10^{01}$ & $ -$ \\ 
				$ 2^3$ & $ 2.1\times10^{00}$& $02.9$ \\
				$ 2^4$ & $ 5.4\times10^{-03}$& $ 08.6$\\
				$ 2^5$ & $ 3.2\times10^{-09}$ & $20.7$ \\
				$ 2^6$ & $ 6.7\times10^{-16}$ & $22.2$  \\ 
				\hline
			\end{tabular}
			\quad \quad
			\begin{tabular}{ c | c | c }
				\hline
				\multicolumn{3}{l}{$\kappa = 200$}\\ \hline
				n & $\varepsilon_{n}$ & noc \\ \hline
				$2^2$ & $ 1.3\times10^{-02}$ & $--$ \\ 
				$2^3$ & $ 6.1\times10^{-03}$ & $01.1$ \\ 
				$2^4$ & $2.2\times10^{-04}$ & $ 04.8$ \\ 
				$2^5$ & $6.0\times10^{-10}$ & $18.5$ \\ 
				$2^6$ & $2.3\times10^{-16}$ & $21.3$ \\ \hline
			\end{tabular}
			\caption{Convergence study for the solution of \cref{eq:yukawa} with density $f$ given in \cref{eq:yukaw}.}\label{table:yukawa}
		\end{center}
	\end{table}
\end{example}

\subsection{Application to the volumetric scattering problem}
~\label{sec:appl-volum-scatt}
In this section, utilizing our integration scheme, we design a fast and accurate Nystr\"{o}m solver for the problem of wave scattering by penetrable inhomogeneous media. 
This solver converges super algebraically when the material properties are smooth. In addition,
our algorithm converges with second-order using the Fourier smoothing technique discussed in \Cref{sec:FS} when discontinuous material interfaces are present. Moreover, this quadratic rate of convergence 
is not constrained by the smoothness of scattering geometry, and the method yields second-order convergence even for obstacles with non-smooth boundaries.
%



Formally, the problem of scattering of time-harmonic acoustic waves by a bounded inhomogeneity $\Omega \subset \mathbb{R}^2$ can be described as follows: for a given incident wave $u^{inc}$ satisfying the free space Helmholtz equation 
\begin{equation}\label{eq:helmholtz}
\Delta u^{inc}(x)+\kappa^2 u^{inc}(x) = 0 \hspace{0.5cm} x\in \mathbb{R}^2,
\end{equation}
find the total field $u$ such that \cite{colton1998inverse}
\begin{equation}\label{eq:totalfield}
\Delta u(x) + \kappa^{2}\mu(x) u(x) = 0, \hspace{0.5cm} x\in \mathbb{R}^2,
\end{equation}
and scattered field $u^{sc} = u-u^{inc}$ satisfying the Sommerfeld radiation condition \cite{colton1998inverse}
\begin{equation}\label{eq:rc}
\lim_{r \to \infty} r^{1/2} \left( \frac{\partial u^{sc}}{\partial r} - i\kappa u^{sc} \right)
= 0,
\end{equation}
where $\kappa$ is the wave number of $u^{inc}$, and $r = |x|$ for $x \in \mathbb{R}^2$. Here, $\mu(x)$  denotes the refractive index of the inhomogeneity that is assumed to be smooth within $\Omega$ and takes the value one outside $\Omega$. Moreover, $\mu$ is allowed to have a jump discontinuity across the material interface $\partial \Omega$.

An equivalent integral equation formulation of
the scattering problem \cref{eq:helmholtz,eq:totalfield,eq:rc} is given by the Lippmann-Schwinger equation \cite{colton1998inverse, martin2003acoustic}
\begin{equation}\label{eq:lippsch}
u(x) + \kappa^{2} (A(mu))(x)= u^{inc}(x),  \quad  x\in\Omega,
\end{equation}
where 
\begin{equation}\label{eq:lippsch-ieo}
(Av)(x) = \int_\Omega  g_\kappa(x-y)v(y) dy
\end{equation}
with the contrast function $m(x) = 1-\mu(x)$ and 
$g_\kappa(x) = \frac{i}{4}H_0^{(1)}(\kappa |x|)$. Without loss of generality, we  take $\Omega \subseteq D$.

An application of \cref{eq:An} at the Nystr\"{o}m nodes  $ x_j = j/n, \, \, j\in \mathbb{G}_n $ for approximating $A(mu)$ in \cref{eq:lippsch} yields the following linear system
\begin{equation}\label{eq:descretls}
u(x_j)+ \kappa^{2} (A_n(mu))(x_j) = u^{inc}(x_j), \quad j\in \mathbb{G}_n.
\end{equation}
Finally, the solution of the linear system \cref{eq:descretls} is obtained using the matrix-free implementation of GMRES in $O(N\log N)$ cost per iterations. 
In the next few examples, we present a variety of scattering calculations using our solver.  
In all the reported experiments, the incident wave is taken as  $u^{inc}(x) = e^{i\kappa x_1}$.

\begin{example}\label{Exp:LippmannOp} \emph{(Scattering by a disc with smooth material properties)}
	As a first exercise of this section, we compute numerical solution of the Lippmann-Schwinger equation \cref{eq:lippsch} over inhomogeneous circular region with acoustical size $\kappa a = 1$, 
	where $a$ denotes the diameter of the inhomogeneity $\Omega$. 
	The refractive index of the inhomogeneous region is given by
	\begin{equation} \label{eq::refLipp}
	\mu(x) =
	\begin{cases}
	2, & |x-x_c| \le t_0, \\
	1+\varkappa(\frac{|x-x_c|-t_{0}}{t_{1}-t_{0}}), & t_0 < |x-x_c| < t_1, \\
	1, & \text{otherwise},
	\end{cases}
	\end{equation}
	where
	\begin{equation} \label{eq::halfwindow1}
	\varkappa(t) =
	\text{exp}\left(\frac{2e^{-1/t}}{t-1}\right)
	\end{equation}
	with $t_1=0.45$, $t_0=0.9 t_1$ and $x_c = (1/2,1/2)$. For a pictorial visualization of the contrast function, we display a plot of $1-\mu(x)$ in \cref{fig:disc_exp1}(a). To study the convergence, numerical solution of  \cref{eq:lippsch} is computed for several level of discretization and the errors are reported in \Cref{table:SolverLipp1}.  The errors from a similar experiment where $\kappa a = 50$ is also reported in \Cref{table:SolverLipp1}.
	The results clearly show the high-order character of our algorithm when  the contrast function $1-\mu$ is smooth and compactly supported.
	A plot of the absolute value of the total field for a disc of acoustical size $\kappa a = 80$ is displayed in \cref{fig:disc_exp1}(b).
	
	\begin{figure}[h!] 
		\centering
		\subfloat[Contrast function] {
			\includegraphics[scale=0.37]{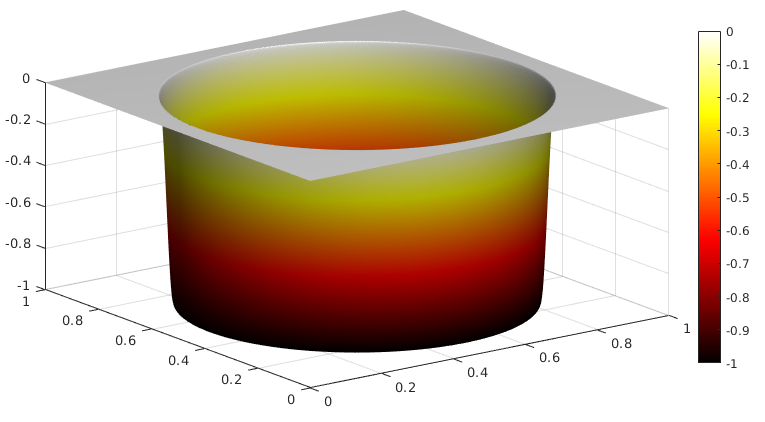}
		}\label{fig:exp1a}
		\subfloat[$|u|$] {
			\includegraphics[scale=0.35]{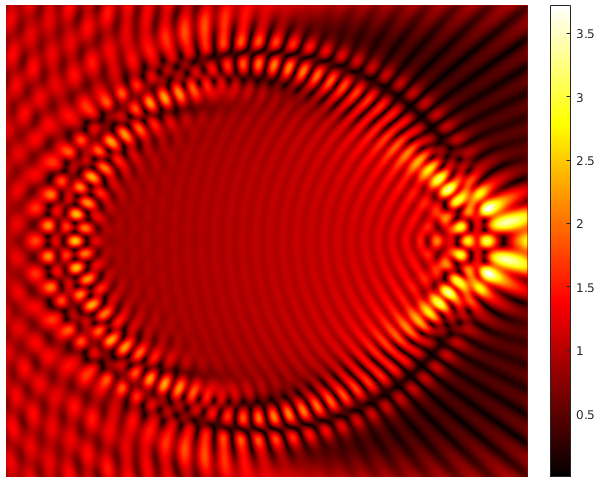}
		}\label{fig:exp1b}
		\caption{(left) The contrast function $1-\mu(x)$ with refractive index $\mu$ given in equation \cref{eq::refLipp}. (right) Solution to the scattering problem in \Cref{Exp:LippmannOp} for the circular disc of size $\kappa a = 80$.}\label{fig:disc_exp1}
	\end{figure}
	
	\begin{table}[t]
		\begin{center}
			\begin{tabular}{ c | c | c }
				\hline
				\multicolumn{3}{l}{$\kappa a = 1$}\\ \hline
				n & $\varepsilon_{n}$ & noc \\ \hline
				$2^5$ & $ 2.0\times10^{-4}$ & $--$ \\ \hline
				$2^6$ & $ 4.1\times10^{-5}$ & $2.3$ \\ \hline
				$2^7$ & $7.5\times10^{-7}$ & $ 5.8$ \\ \hline
				$2^8$ & $1.2\times10^{-8}$ & $6.0$ \\ \hline
				$2^9$ & $4.6\times10^{-11}$ & $8.1$ \\ \hline
			\end{tabular}
			\quad \quad
			\begin{tabular}{ c | c | c }
				\hline
				\multicolumn{3}{l}{$\kappa a= 50$}\\ \hline
				n & $\varepsilon_{n}$ & noc \\ \hline
				$2^5$ & $ 3.5\times10^{-1}$ & $--$ \\ \hline
				$2^6$ & $ 2.5\times10^{-2}$ & $3.8$ \\ \hline
				$2^7$ & $3.8\times10^{-4}$ & $ 5.9$ \\ \hline 
				$2^8$ & $4.0\times10^{-6}$ & $6.5$ \\ \hline
				$2^9$ & $1.2\times10^{-8}$ & $8.4$ \\ \hline
			\end{tabular}
			\caption{Convergence study: plane wave scattering by a circular inhomogeneity with refractive index given in \cref{eq::refLipp}. }\label{table:SolverLipp1}
		\end{center}
	\end{table} 
	
	%
	%
	%
	%
\end{example}

\begin{figure}[h!] 
	\centering
	\subfloat[Contrast function for smooth filter disc] {
		\includegraphics[scale=0.34]{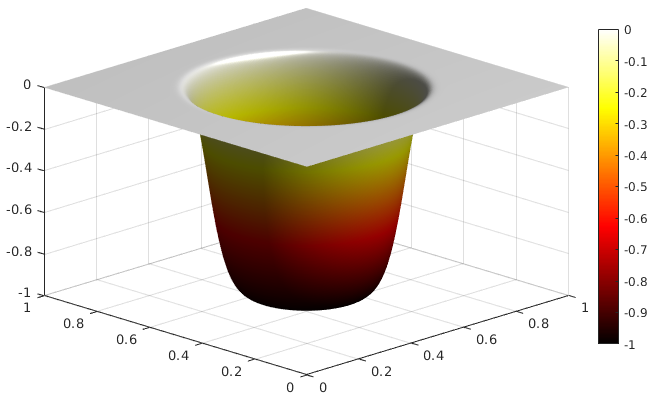}
	}\label{fig:filterD}
	\subfloat[Contrast function for Luneburg lense] {
		\includegraphics[scale=0.32]{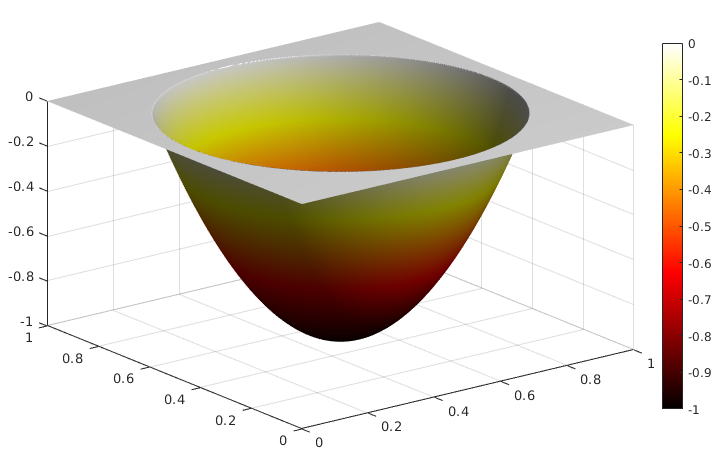}
	}\label{fig:luneberg}
	\caption{ A three dimensional view of the contrast function $m(x)$ for the smooth filter disc and the Luneburg lens defined in \cref{eq:smoothdisc} and \cref{eq:lune}, respectively.
	}\label{fig:contrastf}
\end{figure}

\begin{figure}[h!] 
	\centering
	\subfloat[Scattering by smooth filter disc] {
		\includegraphics[scale=0.14]{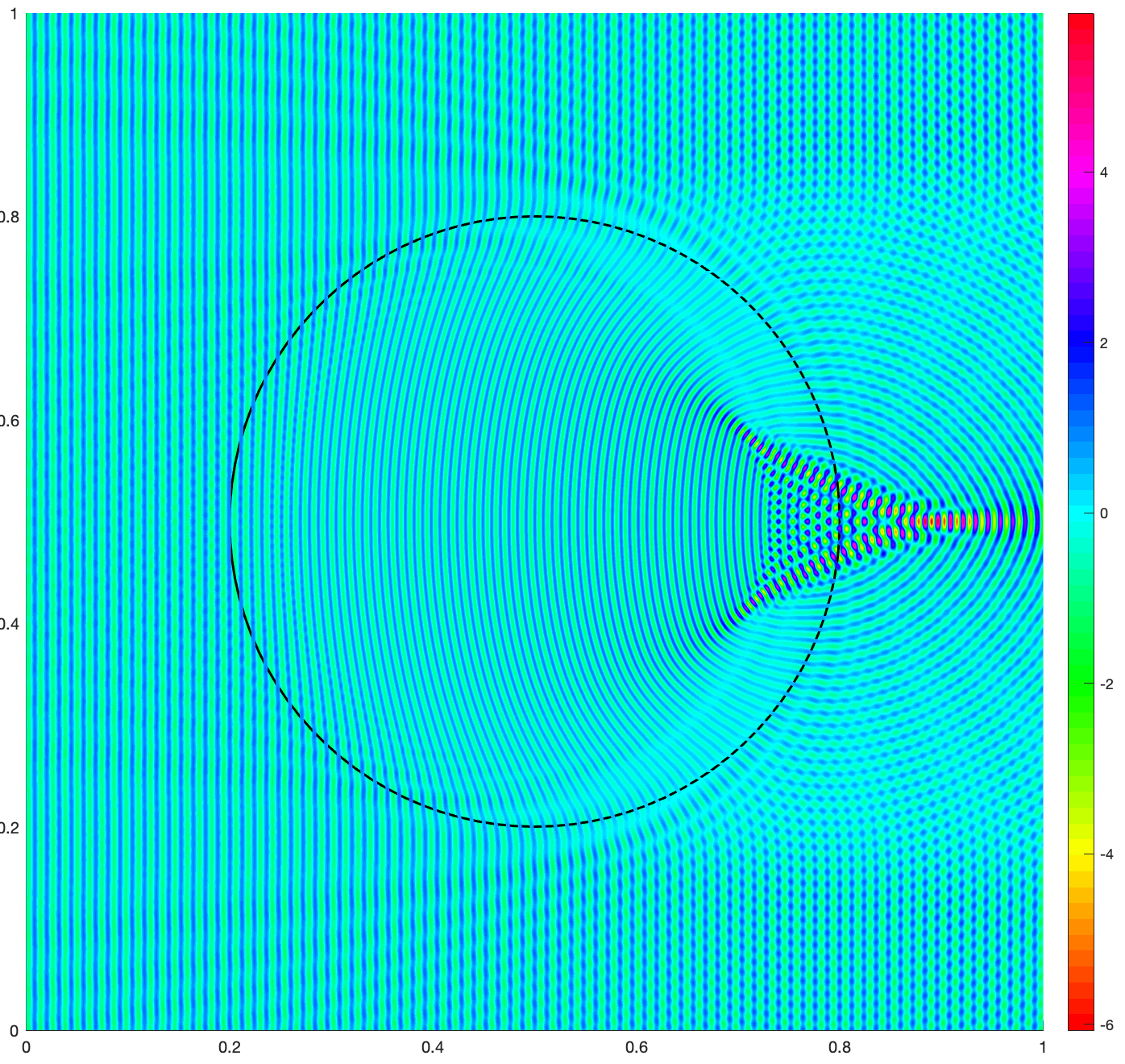}
	}\label{fig:smootha}
	\subfloat[Scattering by Luneburg lense] {
		\includegraphics[scale=0.145]{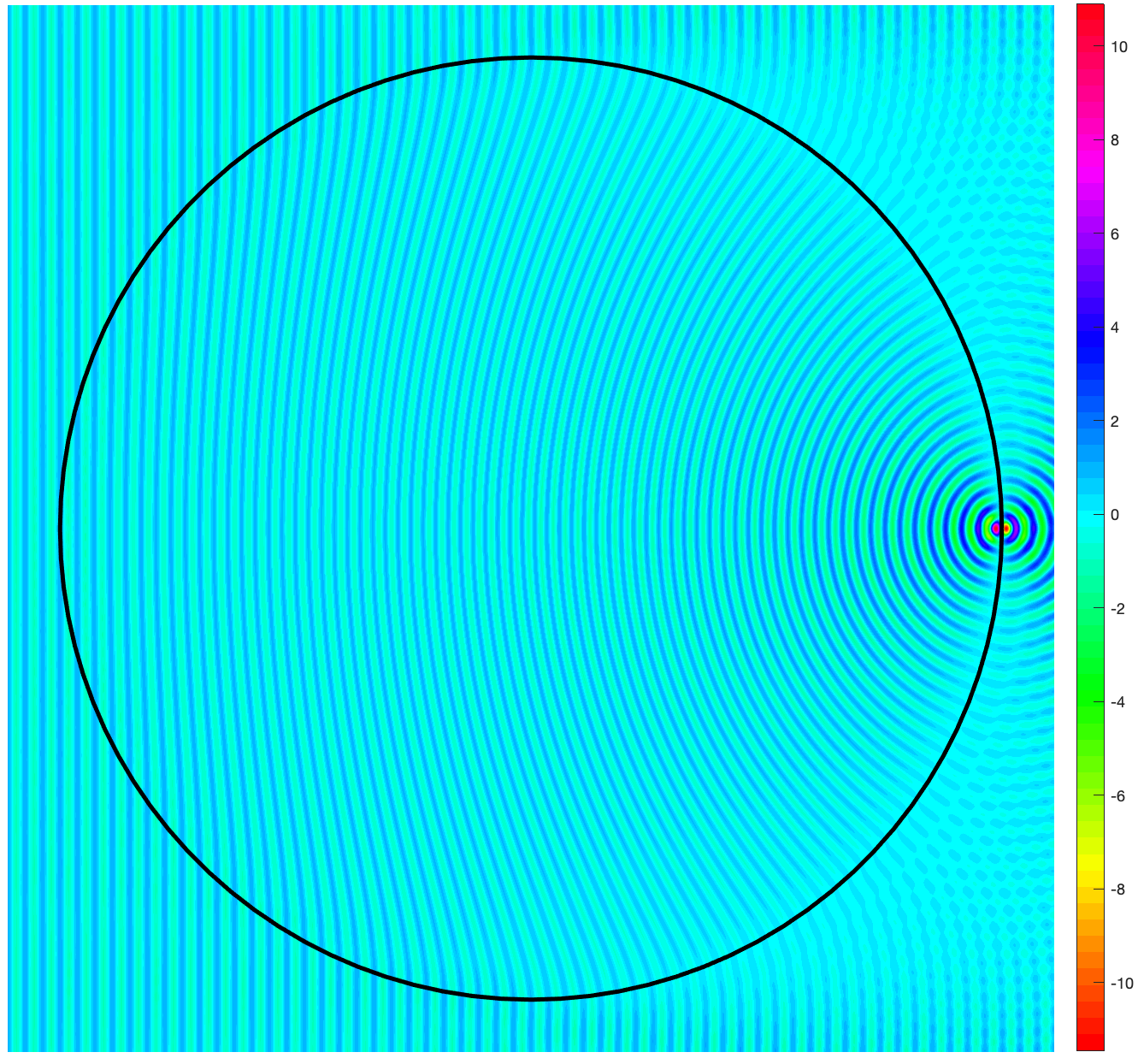}
	}\label{fig:smoothb}
	\caption{(left) Scattering of a plane wave by penetrable smooth filter disc of acoustic size $\kappa a = 80\pi$ in a unit square of size $160\pi$. (right) Scattering of a plane wave by penetrable Luneburg lens of acoustic size $\kappa a = 108\pi$ in a unit square of size $120\pi$. A computational grid of size $n\times n$ with $n=4096$ is used.}\label{fig:smooth1}
\end{figure}

\begin{example}\label{Exp:filtDisc} 
	\emph{(Scattering by a smooth filter disc)}
	In this example, we consider a filter flat disk of diameter $a=0.5$ with center $x_c = (1/2,1/2)$ for which the $m(x)$ is given by \cite{vico2016fast}
	\begin{equation}\label{eq:smoothdisc}
	m(x) = -\exp\left(\frac{-1}{2}\left(\frac{2|x-x_c|}{a}\right)^8\right).
	\end{equation}
	A surface visualisation of the contrast function is shown in \cref{fig:contrastf}(a).
	In \cref{table:smoothdisc}, we show the errors and the corresponding numerical order of convergence for $\kappa a = 2\pi, \, 40\pi, \, 80\pi$.
	The results again show the super algebraic convergence.
	For a pictorial visualization, we display a surface plot of the real part of the computed total field for $\kappa a = 80\pi$ 
	in \cref{fig:smooth1}(a). For this simulation, we use $n = 4096$.
	%
	
	\begin{table}[t!]
		\begin{center}
			\begin{tabular}{ c | c | c |c |c |c |c }
				\hline
				& \multicolumn{2}{l|}{$\kappa a=2\pi$}& \multicolumn{2}{l|}{$\kappa a =40\pi$} & \multicolumn{2}{l}{$\kappa a=80\pi$}\\ \hline
				$n$ & $\varepsilon_{n}$ & noc& $\varepsilon_{n}$ & noc& $\varepsilon_{n}$ & noc \\ \hline
				$2^4$ & $ 4.7\times10^{-4}$ & $-$ & 
				$ -$ & $-$& 
				$ -$ & $-$\\ 
				$2^5$ & $ 7.6\times10^{-7}$ & $9.3$ & 
				$ 1.0\times10^{0}$ & $-$& 
				$ -$ & $-$\\ 
				$2^6$ & $ 2.6\times10^{-11}$ & $14.8$ & 
				$ 1.4\times10^{-2}$ & $6.2$& 
				$ 1.2\times10^{0}$ & $-$\\ 
				$2^7$ & $ 2.6\times10^{-14}$ & $10.0$ & 
				$ 7.9\times10^{-12}$ & $30.7$& 
				$ 2.8\times10^{-4}$ & $12.0$\\ 
				$2^8$ & $ -$ & $-$ & 
				$ 8.6\times10^{-13}$ & $3.2$& 
				$ 1.6\times10^{-12}$ & $27.4$\\ 
				\hline
			\end{tabular}
			\caption{A convergence study for the scattering by smooth filter disc in \Cref{Exp:filtDisc}.}\label{table:smoothdisc}
		\end{center}
	\end{table}	
\end{example}

\begin{example}\label{Exp:lunelens} \emph{(Scattering by Luneburg lens)}
	Next, we consider the nonsmooth Luneburg lens of diameter $a=0.9$, which is designed to focus an incoming wave to a single point. The $m(x)$ for the Luneburg lens is given by \cite{vico2016fast}
	\begin{equation}\label{eq:lune}
	m(x) = \left(\frac{2|x-x_c|}{a} \right)^2 -1,
	\end{equation}
	where $x_c = (1/2,1/2)$. The corresponding contrast function is shown in \cref{fig:contrastf}(b). In \cref{table:lunelens}, we present the computational errors for $\kappa a = 2\pi, \, 40\pi, \, 80\pi$.  
	We display the real part of the total field computed on the $4096\times 4096$ grid for  $\kappa a = 108\pi$ in \cref{fig:smooth1}(b).
	
	\begin{table}[t!]
		\begin{center}
			\begin{tabular}{ c | c | c |c |c |c |c }
				\hline
				& \multicolumn{2}{l|}{$\kappa a=2\pi$}& \multicolumn{2}{l|}{$\kappa a =40\pi$ } & \multicolumn{2}{l}{$\kappa a=80\pi$}\\ \hline
				$n$ & $\varepsilon_{n}$ & noc& $\varepsilon_{n}$ & noc& $\varepsilon_{n}$ & noc \\ \hline
				$2^6$ & $ 4.8\times10^{-5}$ & $-$ & 
				$ 6.5\times10^{-3}$ & $-$& 
				$ 9.3\times10^{-1}$ & $-$\\ 
				$2^7$ & $ 2.2\times10^{-5}$ & $1.1$ & 
				$ 1.1\times10^{-3}$ & $2.5$& 
				$ 4.1\times10^{-3}$ & $7.8$\\ 
				$2^8$ & $ 2.9\times10^{-6}$ & $2.9$ & 
				$ 1.3\times10^{-4}$ & $3.1$& 
				$ 3.7\times10^{-4}$ & $3.5$\\ 
				$2^9$ & $ 1.8\times10^{-7}$ & $4.0$ & 
				$ 3.8\times10^{-5}$ & $1.8$& 
				$ 8.7\times10^{-5}$ & $2.1$\\ 
				$2^{10}$ & $6.6\times10^{-8}$ & $1.5$ & 
				$ 4.5\times10^{-6}$ & $3.1$& 
				$ 7.0\times10^{-6}$ & $3.6$\\ 
				\hline
			\end{tabular}
			\caption{A convergence study for the scattering by Luneburg lens in \Cref{Exp:lunelens}.}\label{table:lunelens}
		\end{center}
	\end{table}	
\end{example}

\begin{figure}[h!] 
	\centering
	\subfloat[Square with cavity] {
		\includegraphics[scale=0.495]{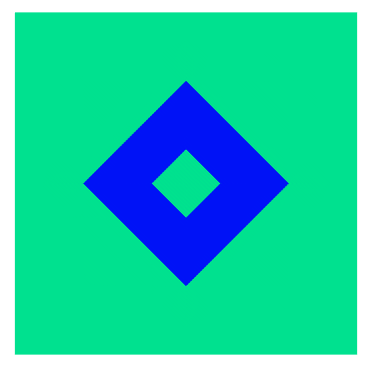}
	} 
	\subfloat[$|u^{sc}|$] {
		\includegraphics[width=0.3\textwidth,trim={3.2cm 1.0cm 2.1cm 1.2cm},clip]{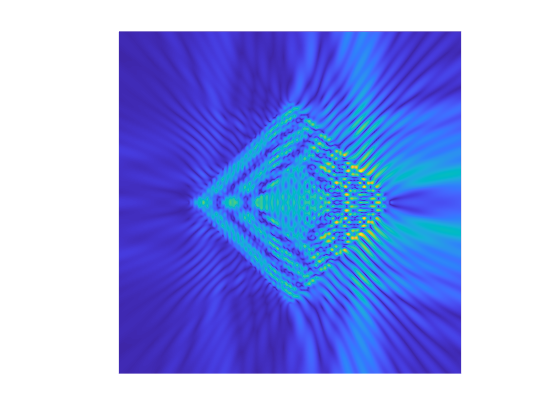}
	}
	\subfloat[$|u|$] {
		\includegraphics[width=0.3\textwidth,trim={3.2cm 1.0cm 2.1cm 1.2cm},clip]{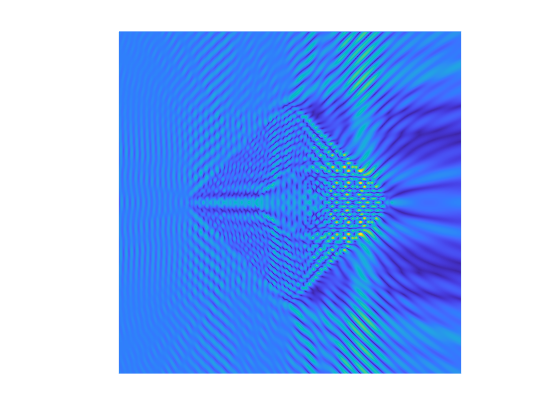}
	}
	\caption{Scattering of a plane wave by penetrable square cavity $\Omega$ of acoustic size $\kappa a = 128$ when contrast function $m(x)=-1$ in $\Omega$ and zero outside. A computational grid of size $512\times 512$ is used for this simulation. }\label{fig:Fsmooth2}
\end{figure}

\begin{example}\label{exp:9}	\emph{(Scattering in the presence of a discontinuous material interface and Fourier smoothing)}
	We now consider some examples where the use of Fourier smoothing has beneficial effect on the rate of convergence of the scattering solver. In particular, through these experiments, we compare the solutions obtained by the Nystr\"om solver with and without Fourier smoothing. 
	
	First, we consider the disc of acoustic size $\kappa a = 40$ with $m(x)=-1$ in the disc and zero outside. A convergence study is shown in \cref{table:Scat_FS}. We observe a significant improvement in the accuracy and enhancement in the  convergence rate for the Fourier smoothed version of the method. 
	
	Finally, to show that the Fourier smoothing technique is not limited to only simple scattering configurations, we consider a  star-shaped  geometry containing cusps and a square cavity with corners.  In \cref{table:Scat_FS} and  \cref{table:Scat_FS1}, we present a numerical convergence study 
	where the acoustic size of the problem is $\kappa a = 60$ and $\kappa a = 64$, respectively. The absolute value of total and scattered fields are displayed in \cref{fig:Fsmooth1} and \cref{fig:Fsmooth2}, where the simulation for the star-shaped geometry uses $\kappa a = 100$ whereas $\kappa a = 128$ in the scattering calculation for the square with a cavity.

	\begin{table}[t!]
		\begin{center}
			\begin{tabular}{ c | c | c |c |c |c |c }
				\hline
				& \multicolumn{4}{l|}{Disc shape scatterer}& \multicolumn{2}{l}{Scatterer with cusp}  \\ 
				\hline
				& \multicolumn{2}{l|}{WFS}& \multicolumn{2}{l|}{FS} & \multicolumn{2}{l}{FS}\\
				\hline
				$n$ & $\varepsilon_{n}$ & noc& $\varepsilon_{n}$& noc & $\varepsilon_{n}$& noc  \\ \hline
				$2^4$ & $ 2.4\times10^{0}$ & $-$ & 
				$ 2.4\times10^{0}$ & $ -$
				& $ -$ & $ -$ \\ 
				$2^5$ & $ 2.1\times10^{0}$ & $0.2$ 
				& $ 7.0\times10^{-1}$& $ 1.8$
				& $ 9.3\times10^{-1}$& $ -$ \\ 
				$2^6$ & $ 5.0\times10^{-1}$ & $ 2.1$ & 
				$ 1.7\times10^{-2}$& $ 5.4$ 
				& $ 6.5\times10^{-2}$& $ 3.9$\\ 
				$2^7$ & $2.0\times10^{-1}$ & $ 1.3$ 
				& $ 3.6\times10^{-3}$ & $2.2$ 
				& $ 5.3\times10^{-3}$& $ 3.6$\\ 
				$2^8$ & $4.7\times10^{-2}$ & $2.0$ & 
				$ 1.1\times10^{-3}$ & $1.7$ 
				& $ 1.2\times10^{-3}$& $ 2.2$ \\ 
				$2^9$ & $1.7\times10^{-2}$ & $1.5$ & 
				$ 2.8\times10^{-4}$ & $2.0$ 
				& $ 2.6\times10^{-4}$& $ 2.2$\\ 
				$2^{10}$ & $4.0\times10^{-3}$ & $2.0$& 
				$ 6.9\times10^{-5}$ &$2.0$ 
				& $ 5.8\times10^{-5}$& $ 2.1$\\ 
				\hline
			\end{tabular}
			\caption{A convergence study for the scattering by a disc with  $\kappa a=40$ and by a geometry with non-smooth boundary of size $ka=60$. The refractive index $\mu(x)$ is discontinuous across the scattering interface.}\label{table:Scat_FS}
		\end{center}
	\end{table}

	\begin{figure}[h!] 
		\centering
		\subfloat[$|u^{sc}|$] {
			\includegraphics[width=0.48\textwidth,trim={3.5cm 1.0cm 2.1cm 1.2cm},clip]{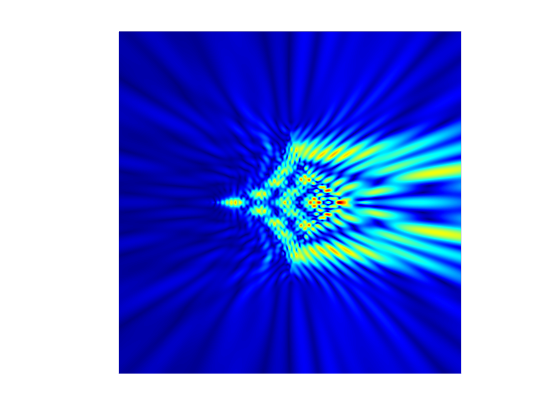}
		}
		\subfloat[$|u|$] {
			\includegraphics[width=0.48\textwidth,trim={3.5cm 1.0cm 2.1cm 1.2cm},clip]{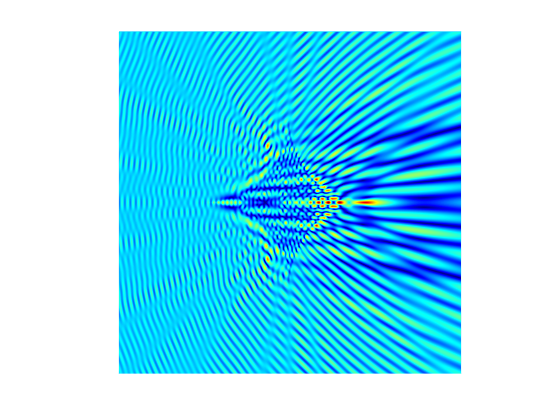}
		}
		\caption{Scattering of a plane wave by penetrable inhomogeneity with cusps of acoustic size $\kappa a = 100$, where $\mu(x)$ is discontinuous across the interface. For this simulation, a computational grid of size $n\times n$ with $n=512$ is used.}\label{fig:Fsmooth1}
	\end{figure}
	
\end{example}

\begin{table}[h] 
	\begin{center}
		\begin{tabular}{ c | c | c |c |c }
			\hline
			& \multicolumn{2}{l|}{WFS}& \multicolumn{2}{l}{FS}\\ \hline
			n & $\varepsilon_{n}$ & noc & $\varepsilon_{n}$ & noc \\ \hline
			$2^5$ & $9.8\times10^{-1}$ & $-$
			& $1.0\times10^{-0}$ & $-$ \\ 
			
			$2^6$ & $4.5\times10^{-1}$ & $1.1$ 
			& $6.2\times10^{-2}$ & $4.0$\\ 
			
			$2^7$ & $1.8\times10^{-1}$ & $1.4$ 
			& $2.9\times10^{-3}$ & $4.4$\\ 
			
			$2^8$ & $1.4\times10^{-1}$ & $0.4$ 
			& $6.8\times10^{-4}$ & $2.0$\\ 
			
			$2^9$ & $7.3\times10^{-2}$ & $0.9$ 
			& $1.6\times10^{-4}$ & $2.1$\\ 
			
			$2^{10}$ & $2.1\times10^{-2}$ & $1.8$ 
			& $4.7\times10^{-5}$ & $1.8$\\\hline
		\end{tabular}
		\caption{A convergence study for a plane wave scattering by penetrable square cavity $\Omega$ in \Cref{exp:9}, where the contrast function $m(x)=-1$ in $\Omega$ and zero outside.}\label{table:Scat_FS1}
	\end{center}
\end{table} 

\section{Concluding discussion and future directions}
~\label{sec:conclusions}
In this article, we present a simple, efficient, and high-order numerical integration scheme to evaluate convolution integrals of the form
\begin{equation*}
(Au)(x) = \displaystyle \int_D g(x-y)u(y)\,dy,
\end{equation*}
where $g$ is weakly singular and  $u$ is supported in $D$. The preeminent motivation for this numerical integration scheme is to compute the convolution $Au$ on a uniform grid of size $N$ in $O(N\log N)$ operations with high-order accuracy. To achieve this, we primarily rely on a periodic Fourier extension $u_e$ of $u$ with a suitably large period $b$. 
We demonstrate the accuracy and efficiency of the proposed integration scheme through a variety of numerical experiments.
We see that the rate of convergence of the method increases with increasing smoothness of the periodic extension. In fact, the method exhibits super-algebraic convergence when the extension is infinitely differentiable. Furthermore, when the density has jump discontinuities, we employ a certain Fourier smoothing technique to accelerate the convergence to achieve the quadratic rate.  
We utilize the integration scheme for numerical solution of certain partial differential equations and to obtain a fast and high-order Nystöm solver for the solution of the Lippmann-Schwinger integral equation.

We observe that the rate of convergence of the proposed scheme depends on $g$ and the smoothness of $u_e$.
The numerical experiments indicate that the computational error behaves as
\begin{align*}
\Vert Au - A_nu \Vert_{\infty} &\le 
O\left(n^{-(2+r+\min{(0, \beta-1)})}\right), 
\end{align*} 
provided $u_e$ is $r$ times continuously differentiable in $\mathbb{R}^2$ with Lipschitz continuous $u_e^{(r+1)}$, that is,
$u_{e}\in C^{r,1}(\mathbb{R}^2)$ and $g$ is such that $|\widehat{g}(k)| \le C_g (1+|k|^2)^{-\beta/2}$ for some constant $C_g$.
The error analysis to rigorously establish this behaviour is an interesting future research effort. The implementation and analysis of the method to three and higher dimensions is another future research direction of significant interest.

\section*{Acknowledgments}
AA acknowledges the support by Science \& Engineering Research Board through File No. MTR/2017/000643.
AKT \& AP acknowledge the Initiation Grant from the Indian Institute of Science Education and Research Bhopal.

\end{document}